\documentclass[11pt]{amsart}

\usepackage[
  margin=1in
]{geometry}

\usepackage{amscd}
\usepackage{amssymb}
\usepackage{amsmath}
\usepackage{amsfonts}
\usepackage{enumerate}
\usepackage[all,cmtip]{xy}
\usepackage{tikz-cd}
\usepackage{hyperref}
\usepackage{xcolor}
\usepackage{adjustbox}
\usepackage{graphicx}
\usepackage[T1]{fontenc}
\usepackage{lmodern}

\newtheorem{thm}{Theorem}[section]
\newtheorem{prop}[thm]{Proposition}
\newtheorem{lem}[thm]{Lemma}
\newtheorem{cor}[thm]{Corollary}

\theoremstyle{definition}

\newcommand{\Ss}{\mathbb{S}}
\newcommand{\Zz}{\mathbb{Z}}
\newcommand{\Ee}{\mathbb{E}}
\newcommand{\Ff}{\overline{\mathbb{F}}_p}
\newcommand{\BP}[1]{\mathrm{BP}\langle #1 \rangle}
\newcommand{\THH}{\mathrm{THH}}
\newcommand{\TT}[1]{T^{#1}\otimes}

\newcommand{\MU}{\mathrm{MU}}
\newcommand{\Wn}{W(n,m)}
\newcommand{\Rr}{\pi_0(R)_p^{\wedge}}
\newcommand{\threerightarrows}{
\mathrel{\vcenter{\hbox{$\rightrightarrows$}\nointerlineskip\hbox{$\rightarrow$}}}
}

\makeatletter
\@namedef{subjclassname@2020}{2020 Mathematics Subject Classification}
\makeatother


\numberwithin{equation}{section}

\title[Height of $(\TT{n}R)^{hT^n}$ with $R$ connective]{Chromatic heights of homotopy fixed points of iterated topological Hochschild homology of connective spectra}

\author[J. Zha]{Jiaxi Zha}
\address{Department of Mathematics, Nankai University, No.94 Weijin Road, Tianjin 300071, P. R. China}
\email{1093913699@qq.com}

\begin{document}

\begin{abstract}
	For a prime $p$, let $R$ be a $p$-local $\Ee_{\infty}$-ring spectrum of height $m$. When $m\geq0$, the chromatic redshift theorem implies that its $n$-fold iterated algebraic $K$-theory $K^{(n)}(R)$ has height $m+n$. Moreover, there is a natural map $K^{(n)}(R)\to (T^n\otimes R)^{hT^n}$, where $T^n\otimes R$ denotes the $n$-fold iterated topological Hochschild homology. Hence it is natural to ask whether the target also has height $m+n$. Although this fails in general, we prove that it holds whenever $R$ is connective and $m=-1$ or $m\geq1$. We also determine the height of $(T^n\otimes R)^{hT^n}$ when $R$ is connective and $m=0$. Our proof is based on an iterative calculation using Wilson spaces and relative topological Hochschild homology. As an application, we determine the support of the $T^n$-homotopy fixed points of $n$-fold iterated topological Hochschild homology of $\BP{m}$ for $n=2,3$ with $\BP{m}$ equipped with any $\Ee_3$-$\MU$-algebra form.
\end{abstract}

\maketitle

\section{Introduction}
Fix a prime $p$ and a finite $p$-local spectrum $U$ of type $m$. Let $K(m)$ be the height $m$ Morava $K$-theory. By \cite[Theorem 9]{HS98}, there is a $v_m$-self map $v\colon \Sigma^dU\to U$ which induces an isomorphism after taking $K(m)$-homology. Set $T(m)$ to be $U[v^{-1}]$. The support of a $p$-local spectrum $X$ is defined to be
\[\mathrm{supp}(X):=\left\{m\in\mathbb{N}\ |\ T(m)\otimes X\not\simeq0\right\}.\]
When $X$ admits an $\Ee_{\infty}$-ring structure, Hahn \cite{Hahn16} proved that the support of $X$ is either empty or an interval containing $0$. Therefore, it is natural to define 
\[\mathrm{height}(X):=\max\left\{m\geq-1\ |\ T(m)\otimes X\not\simeq0\right\},\]
where $T(-1)$ is taken to be $\Ss_{(p)}$. Moreover, \cite[Lemma 2.3]{LMMT24} implies that we can replace $T(m)$ by $K(m)$ in this case.

The chromatic redshift theorem \cite[Theorem E]{BSY22} implies that the algebraic $K$-theory spectrum of an $\Ee_{\infty}$-ring spectrum of height $m\geq0$ has height $m+1$. By \cite[Theorem 7.4]{BGT14}, for an $\Ee_{\infty}$-ring spectrum, there is a natural map of $\Ee_{\infty}$-ring spectra $K(R)\to\THH(R)^{hT}$. Consequently, one might hope that $\THH(-)^{hT}$ also raises height by one for $\Ee_{\infty}$-ring spectra. In general, it fails. For example, $\THH(E_m)^{hT}$ still has height $m$, where $E_m$ is the height $m$ Morava $E$-theory; see \cite[Example 5.2]{Zha26}. However, the computations of Ausoni \cite{Aus10} and Ausoni--Rognes \cite{AR02} show that this does hold for the Adams summand and connective complex $K$-theory for $p\geq5$. Moreover, in a draft on higher redshift, Rognes showed that $p\in\pi_0(\Ss)$, $v_1\in\pi_{2p-2}(\Ss/p)$ and $v_2\in\pi_{2p^2-2}(V(1))$ have nonzero images in $\pi_*((T\otimes\mathbb{F}_p)^{hT})$, $V(0)_*((T^2\otimes\mathbb{F}_p)^{hT^2})$ and $V(1)_*((T^3\otimes\mathbb{F}_p)^{hT^3})$, respectively, with $p$ odd in the third case. Motivated by this draft, Veen \cite{Vee18} subsequently proved that, when $p\geq5$ and $1\leq n\leq p$ or $p=3$ and $1\leq n\leq2$, the unit map $k(n-1)_*\to k(n-1)_*((T^n\otimes\mathbb{F}_p)^{hT^n})$ sends $v_{n-1}$ to a nonzero element, where $k(n-1)$ is the connective height $n-1$ Morava $K$-theory and $T^n\otimes\mathbb{F}_p$ is the $n$-fold iterated topological Hochschild homology; see Section \ref{sec2} for an explanation of this notation. 

In recent work \cite{Zha26}, the author extended these height computations and showed that the heights of $(T^n\otimes\mathbb{F}_p)^{hT^n}$, $(T^n\otimes\Zz_{(p)})^{hT^n}$ and $(T^n\otimes j_p)^{hT^n}$ are $n-1$, $n$ and $n+1$, respectively, for any $n$ and prime $p$. The author's method differs from those used in pervious work. Earlier computations used spectral sequences to detect images of periodic elements. In contrast, the author used relative topological Hochschild homology and the Wilson space method developed by Hahn--Raksit--Wilson in \cite{HRW25} for their construction of the motivic filtration on topological cyclic homology. These calculations lead us to ask whether the functor $(T^n\otimes -)^{hT^n}$ raises the height of a $p$-local connective $\mathbb E_\infty$-ring spectrum $R$ by $n$. Our first main result shows that it does whenever $R$ has nonzero height.

\begin{thm}\label{thm1}
	For $m=-1$ or $m\geq1$, let $R$ be a $p$-local\footnote{Since the proof of Theorem \ref{thm1} does not use the assumption that $R$ is $p$-local when $m\geq 1$, this assumption can be dropped in this case.} connective $\Ee_{\infty}$-ring spectrum of height $m$. For any $n\geq1$,  the height of $(T^n\otimes R)^{hT^n}$ is $m+n$.
\end{thm}

Indeed, tensoring with $K(k)$, the cyclotomic trace map, the canonical natural transformation $\mathrm{TC}(-)\to\THH(-)^{hT}$, and the comparison map of limits give an $\Ee_1$-map 
\[K(k)\otimes K^{(n)}(R)\longrightarrow K(k)\otimes(T^n\otimes R)^{hT^n}\]
by \cite[Theorem 7.4]{BGT14}, where $K^{(n)}(R)$ is the $n$-fold iterated algebraic $K$-theory of $R$. By \cite[Theorem E]{BSY22}, $K^{(n)}(R)$ has height $m+n$. Hence it follows from \cite[Theorem 1.1]{Hahn16} that we only need to show that $K(m+n)\otimes(T^n\otimes R)^{hT^n}$ is nonzero. 

For $m=-1$, we reduce this to the computation of the height of $(T^n\otimes\mathbb{F}_p)^{hT^n}$ in \cite{Zha26}.  For $m\geq1$, the chromatic Nullstellensatz of Burklund--Schlank--Yuan \cite{BSY22} gives a map of $\Ee_{\infty}$-ring spectra $R\to\tau_{\geq0}E_m(L)$. We prove this by showing that $(T^n\otimes\tau_{\geq0}E_m(L))^{hT^n}$ has height $m+n$, where $L$ is an algebraically closed field of characteristic $p$. The case $m=0$ is different, since the height of $(T^n\otimes R)^{hT^n}$ can vary with $R$. For example, take $R$ to be the height zero connective Morava $E$-theory $\tau_{\geq0}E_0(L)$, where $L$ is an algebraically closed field of characteristic zero. Then $(T^n\otimes R)^{hT^n}$ is still rational, since rationalization is a smashing localization. Hence it still has height zero. On the other hand, $(T^n\otimes H\Zz)^{hT^n}$ has height $n$ by \cite[Theorem 1.2]{Zha26}. Nevertheless, we can completely determine the height of $(T^n\otimes R)^{hT^n}$ when $R$ is a $p$-local connective $\Ee_{\infty}$-ring spectrum of height zero. This is our second main theorem. 

\begin{thm}\label{thm2}
	Let $R$ be a $p$-local connective $\Ee_{\infty}$-ring spectrum of height zero. For any $n\geq1$, we have that
	\begin{equation*}
		\mathrm{height}((T^n\otimes R)^{hT^n})=
		\begin{cases}
			0, & p\in\pi_0(R)^{\times}, \\
			n-1, & \text{$p\notin\pi_0(R)^{\times}$ and $\Rr[1/p]=0$,} \\
			n, & \Rr[1/p]\neq0.
		\end{cases}
	\end{equation*}
\end{thm}

To the best of our knowledge, the chromatic redshift theorem is currently known to hold only for $\Ee_{\infty}$-ring spectra\footnote{For example, let $R$ be a connective $\Ee_{\infty}$-ring spectrum of height $m$. Then $\mathrm{End}_R(R^{\oplus\mathbb{N}})$ is a connective $\Ee_1$-ring spectrum of height $m$. However, we have that $K(\mathrm{End}_R(R^{\oplus\mathbb{N}}))\simeq0$ by the Eilenberg swindle.}. However, Hahn--Wilson \cite{HW22} showed that, for a fixed $\Ee_3$-$\MU$-algebra form of $\BP{n}$, both $K(n+1)\otimes K(\BP{n})$ and $K(n+1)\otimes\THH(\BP{n})^{hT}$ are nonzero. Together with the preceding discussion, this motivates the question of whether both $K^{(m)}(\BP{n})$ and $(\int_{T^m}\BP{n})^{hT^m}$ have support $\{0,1,\cdots,m+n\}$ for $m=2,3$; see Section \ref{sec2} for the latter notation. Using Theorem \ref{thm1}, we show that the support of $(\int_{T^m}\BP{n})^{hT^m}$ contains $\{0,1,\cdots,m+n\}$, and we can prove the reverse inclusion holds as well, by the purity theorem of Land--Mathew--Meier--Tamme \cite{LMMT24}. This gives our third main theorem.

\begin{thm}\label{thm3}
	For $n\geq1$, fix any $\Ee_3$-$\MU$-algebra form of $\BP{n}$. The supports of $(\int_{T^2}\BP{n})^{hT^2}$ and $(\int_{T^3}\BP{n})^{hT^3}$ are $\{0,1,\cdots,n+2\}$ and $\{0,1,\cdots,n+3\}$, respectively.
\end{thm}

Finally, we point out that this theorem, together with the purity theorem in \cite{LMMT24}, allows us to determine the support of $K^{(m)}(\BP{n})$ for $m=2,3$; see Corollary \ref{cor:supp(K(BPn))}.

\subsection*{Organization}
In Section \ref{sec2}, we recall some properties of iterated topological Hochschild homology and Morava $E$-theory. In Section \ref{sec3}, we construct an $\Ee_{\infty}$-ring spectrum $\Wn$ and a map of $\Ee_{\infty}$-ring spectra $(T^n\otimes\tau_{\geq0}E_m(L))^{hT^n}\to\Wn^{hT^n}$ by induction on $n$, for $m\geq1$. Then we show that $K(m+n)\otimes\Wn^{hT^n}$ is nonzero and prove Theorem \ref{thm1} in Section \ref{sec4}. Theorem \ref{thm2} is proved in Section \ref{sec5}. Finally, in Section \ref{sec6}, we compute the supports of $(\int_{T^m}\BP{n})^{hT^m}$ and $K^{(m)}(\BP{n})$ for $m=2,3$.

\subsection*{Conventions}
We work with $\infty$-categories throughout, following \cite{HTT,HA}. Therefore, all categories are understood to be $\infty$-categories unless otherwise specified. Moreover,
\begin{enumerate}
\item We write $\mathcal{S}$ for the category of spaces, and $\mathrm{Sp}$ for the category of spectra. The $\Ee_k$-ring spectra are called $\Ee_k$-rings. The maps of $\Ee_k$-ring are called $\Ee_k$-maps. 
\item For a commutative ring $R$, the Eilenberg-MacLane spectrum $HR$ is also denoted by $R$.
\item $\MU$ denotes the $p$-local complex cobordism spectrum. For a spectrum $X$, we denote $X_p$ by its $p$-completion. The ordinary $p$-adic completion of a commutative ring $A$ is denoted by $A_p^{\wedge}$.
\item We fix an $\Ee_3$-$\MU$-form of $\BP{n}$ for each $n\geq1$, which exists by \cite[Theorem 2.0.6]{HW22}. 
\item We write $\Lambda_R(x)$ and $\Gamma_R(x)$ for the exterior and divided power algebras over a commutative ring $R$, respectively.
\item In a commutative ring $R$ and $x,y\in R$, we write $x\doteq y$ to mean that there is a unit $u\in R^{\times}$ such that $x=uy$.
\end{enumerate}

\subsection*{Acknowledgements}
The author is grateful to his advisor Xiangjun Wang for his guidance and support. The author also thanks John Rognes for sharing his draft on higher redshift and for providing historical background for this work.

\subsection*{AI Statement.}
Building on the methods of the author's earlier work \cite{Zha26}, the author conceived the main arguments and developed the original proofs. One exception is the proof of Lemma \ref{lem:A to B}. In the proof of Theorem \ref{thm2}, the author proposed to find a suitable map to reduce the problem to the case of an Eilenberg--MacLane spectrum. Following this idea, GPT-6.2 Astra formulated Lemma \ref{lem:A to B} and supplied its proof. The author also used GPT-6.1 Sol to assist with language editing, identify typographical errors, find potential gaps in the proofs, and search relevant references. The author reviewed the AI-generated suggestions and takes full responsibility for the content of the manuscript.

\section{Recollections on iterated THH and Morava $E$-theory}\label{sec2}
In this section, we recall some properties of iterated THH and Morava $E$-theory which will be used later. Let $\mathcal{C}$ be a presentable symmetric monoidal category. We begin with equivalent descriptions of iterated THH. Let $\mathrm{Mfld}^{\mathrm{fr}}_n$ be the category of framed $n$-manifolds as constructed in \cite{AF15}. There is a full subcategory of $\mathrm{Mfld}^{\mathrm{fr}}_n$ spanned by $n$-manifolds equivalent to finite disjoint unions of copies of $\mathbb{R}^n$, denoted by $\mathrm{Disk}_n^{\mathrm{fr}}$. The factorization homology is defined by the left Kan extension along the inclusion $\mathrm{Disk}_n^{\mathrm{fr}}\hookrightarrow\mathrm{Mfld}^{\mathrm{fr}}_n$
\[\int\colon\mathrm{Alg}_{\Ee_n}(\mathcal{C})\simeq\mathrm{Fun}^{\otimes}(\mathrm{Disk}_n^{\mathrm{fr}},\mathcal{C})\longrightarrow\mathrm{Fun}^{\otimes}(\mathrm{Mfld}^{\mathrm{fr}}_n,\mathcal{C}),\]
where $\mathrm{Fun}^{\otimes}(-,-)$ denotes symmetric monoidal functors. For an $\Ee_n$-algebra $A$ in $\mathcal{C}$, \cite[Lemma 3.18]{AF15} implies that there is an equivalence 
\[\int_{T^n}A\simeq\int_{\mathbb{R}\times T^{n-1}}A\otimes_{\int_{(\mathbb{R}\cup\mathbb{R})\times T^{n-1}}A}\int_{\mathbb{R}\times T^{n-1}}A\simeq \THH(\int_{\mathbb{R}\times T^{n-1}}A).\]
In particular, $\int_{T^n}A$ is the $n$-fold iterated THH of $A$. We will apply this description to the iterated THH for $\BP{n}$. For $\Ee_{\infty}$-algebras, iterated THH can also be described via tensoring with spaces. By \cite[Corollary 3.2.3.5]{HA}, $\mathrm{CAlg}(\mathcal{C})$ is presentable. Hence \cite[Corollary 4.4.4.9]{HTT} gives a functor 
\[\otimes_{\mathcal{C}}\colon \mathcal{S}\times\mathrm{CAlg}(\mathcal{C})\longrightarrow \mathrm{CAlg}(\mathcal{C}),\]
which sends a pair $(X,A)$ to $\mathrm{colim}A_X$. Here $A_X\colon X\to\mathrm{CAlg}(\mathcal{C})$ is the constant functor with the value $A$. For a fixed $\Ee_{\infty}$-algebra $A$ in $\mathcal{C}$, the functor $-\otimes_{\mathcal{C}}A$ preserves colimits by \cite[Example 2.5]{HY17}. Consequently, we also obtain that $\TT{n}A\simeq\THH(\TT{n-1}A).$ We will apply this description to the iterated THH for $\Ee_{\infty}$-rings and $\Ee_{\infty}$-$\MU$-algebras. For an $\Ee_{\infty}$-ring $R$, as the right adjoint of the symmetric monoidal functor $-\otimes R\colon \mathrm{Sp}\to\mathrm{Mod}_R$, the forgetful functor is lax symmetric monoidal. Consequently, it induces a functor $\mathrm{CAlg}(\mathrm{Mod}_R)\to\mathrm{CAlg}(\mathrm{Sp})$. In particular, for $X\in\mathcal{S}$ and $A\in\mathrm{CAlg}(\mathrm{Mod}_R)$, we obtain a comparison map
\begin{equation}\label{eq:forget}
	X\otimes A\to X\otimes_RA,
\end{equation}
of $\Ee_{\infty}$-rings, where $A$ on the left is considered as an $\Ee_{\infty}$-ring via the forgetful functor.

Next, we recall some properties of Morava $E$-theory. From now on, we fix a prime $p$ and $m\geq1$. Given a perfect field $k$ of characteristic $p$ and a formal group $\mathbb{G}$ of height $m$ over $k$, there is an $\Ee_{\infty}$-ring $E(k;\mathbb{G})$ whose homotopy group is isomorphic to
\[\pi_*(E(k;\mathbb{G}))\cong W(k)[[u_1,\cdots,u_{m-1}]][u^{\pm1}],\]
where $W(k)$ is the $p$-typical Witt vectors of $k$, $u_1,\cdots,u_{m-1}$ have degree zero, and $u$ has degree two; see \cite[Definition 2.32]{BSY22} and \cite[\S5]{Lur18}. For an algebraically closed field $L$ of characteristic $p$, we simply write $E_m(L)$ for $E(L;\mathbb{G})$, since formal groups of height $m$ over an algebraically closed field are isomorphic to each other.

Now, we make a reduction of Theorem \ref{thm1}. Let $e_m(L)$ denote $\tau_{\geq0}E_m(L)$ for an algebraically closed field $L$ of characteristic $p$. Lemma 2.2 and Lemma 2.3 in \cite{LMMT24} imply that the height of $e_m(L)$ is the same as that of $E_m(L)$, which is $m$ by \cite[Theorem 5.0.2]{Lur18}. For a connective $\Ee_{\infty}$-ring $R$ of height $m$,  $L_{T(m)}R$ is nonzero. Hence \cite[Theorem A]{BSY22} gives an $\Ee_{\infty}$-map 
\[R\to L_{T(m)}R\to E_m(L)\]
for some algebraically closed field $L$, which further gives an $\Ee_1$-map
\[K(i)\otimes(\TT{n}R)^{hT^n}\longrightarrow K(i)\otimes(\TT{n}e_m(L))^{hT^n},\]
for $i,n\geq1$. Therefore, we deduce that
\[\mathrm{height}((\TT{n}R)^{hT^n})\geq\mathrm{height}((\TT{n}e_m(L))^{hT^n}).\]
As observed in the introduction, the chromatic redshift implies that $\mathrm{height}((\TT{n}R)^{hT^n})\leq m+n$. Consequently, to prove Theorem \ref{thm1}, it suffices to show that 
\[K(m+n)\otimes(\TT{n}e_m(L))^{hT^n}\neq0\]
by the discussion in the introduction again.

We end this section with the following lemma which will be used in the proof of Theorem \ref{thm3}.

\begin{lem}\label{lem:e3}
	Let $R$ be an $\Ee_3$-ring whose homotopy group is isomorphic to $W(\Ff)[[u_1,\cdots,u_{m-1}]][u^{\pm1}]$. If the corresponding formal group law over $\pi_0(R)$ is a universal deformation of a height $m$ formal group law over $\Ff$, then $R$ is equivalent to the underlying $\Ee_3$-ring of $E_m(\Ff)$.
\end{lem}
\begin{proof}
	Since all height $m$ formal groups over $\Ff$ are isomorphic, the underlying $\Ee_1$-ring of $R$ is equivalent to the underlying $\Ee_1$-ring of $E_m(\Ff)$ by \cite[Theorem 7.1]{Rezk98}. In particular, their images in $\mathrm{Alg}(\mathrm{ho}(\mathrm{Sp}))$ are the same. Hence \cite[Corollary 4.53]{Ram23} gives the desired result.
\end{proof}

\section{The $\Ee_{\infty}$-ring $\Wn$}\label{sec3}
In this section, we construct the $\Ee_{\infty}$-ring $W(n,m)$ and establish some of its properties. The construction is similar to that in \cite[\S3]{Zha26}. First, we need some $\Ee_{\infty}$-maps. Let $\mathbb{W}(L)$ be the $p$-complete spherical Witt vectors of an algebraically closed field $L$. By the universal property of the spherical Witt vectors there is an equivalence
\[\mathrm{Map}_{\mathrm{CAlg}(\mathrm{Sp}_p)}(\mathbb{W}(L),e_m(L))\simeq\mathrm{Map}_{\mathrm{CRing}}(L,L[[u_1,\cdots,u_{m-1}]]),\]
see \cite[Example 5.2.7]{Lur18}. Hence we obtain an $\Ee_{\infty}$-map $\mathbb{W}(L)\to e_m(L)$ corresponding to the inclusion $L\hookrightarrow L[[u_1,\cdots,u_{m-1}]]$. Base change along the unit map $\Ss\to\MU$, we obtain the first map of $\Ee_{\infty}$-$\MU$-algebras
\[\MU\otimes\mathbb{W}(L)\longrightarrow e_m(L).\]
The second map is the $\Ee_{\infty}$-$\MU$ orientation
\begin{equation}\label{eq:orient}
	\MU\longrightarrow e_m(L),
\end{equation}
as in \cite[Corollary 8.13]{BSY22}. By \cite[\S6.9]{Rezk98}, we can choose another regular sequence $w_1,\cdots,w_{m-1}$ in $\pi_0(E_m(L))$ such that $v_i\in\MU_*$ is sent to $w_iu^{p^i-1}$ for $1\leq i\leq m-1$. Moreover, there is an isomorphism $\pi_*(e_m(L))\cong W(L)[[w_1,\cdots,w_{m-1}]][u]$. Denote the image of $v_i$ by $c_iu^{p^i-1}$ for $i\geq m$. Since the height of the formal group law is $m$, $c_m$ is a unit in $\pi_0(e_m(L))$. By abuse of notation, we still use the notation $u_1,\ldots,u_{m-1}$, with $v_i$ sent to $u_i u^{p^i-1}$ for $1\leq i\leq m-1$. 

Moreover, let $W_{2k}$ be the Wilson space as in \cite[Definition 3.1.1]{HRW25} for $k\geq0$. Since $u_i+1$ is a unit in $\pi_0(e_m(L))$, base change along the map \eqref{eq:orient} and using \cite[Proposition 3.1.5]{HRW25} imply that there is a map of $\Ee_{\infty}$-$\MU$-algebras
\[\MU\otimes\Sigma^{\infty}_+W_0\longrightarrow e_m(L),\]
which sends the generator of lowest degree $x_i\in\MU_0(W_0)\cong\Zz_{(p)}[x_i^{\pm}]$ to $u_i+1$ for $1\leq i\leq m-1$. The same proposition also implies that there is a map of $\Ee_{\infty}$-$\MU$-algebras
\[\MU\otimes\Sigma^{\infty}_+W_2\longrightarrow e_m(L),\]
which sends $u\in\MU_2(W_2)$ to the same named element in $\pi_*(e_m(L))$.

Multiplying these maps together in the category of $\Ee_{\infty}$-$\MU$-algebras, we obtain a map of $\Ee_{\infty}$-$\MU$-algebras
\[\MU\otimes\mathbb{W}(L)\otimes\Sigma^{\infty}_+(W_0^{\times m-1}\times W_2)\longrightarrow e_m(L).\]
Since $e_m(L)$ is complete with respect to $(p,u_1,\cdots,u_{m-1})$, this map factors through
\begin{equation}\label{eq:W(0,m)}
	f_0\colon W(0,m):=(\MU\otimes\mathbb{W}(L)\otimes\Sigma^{\infty}_+(W_0^{\times m-1}\times W_2))^{\wedge}_{(p,x_1-1,\cdots,x_{m-1}-1)}\to e_m(L).
\end{equation}

\begin{lem}\label{lem:pi(W0)}
	There is an isomorphism of commutative rings
	\[\pi_*(W(0,m))\cong W(L)[[u_1,\cdots,u_{m-1}]][u,v_i,y_{\delta}].\]
	Here the classes $y_\delta$ arise from the generators of $\pi_*(MU)$ other than the $v_i$, and from the generators of $MU_*(W_{2k})$ other than the generator of lowest degree and $\pi_*(MU)$; see \cite[Theorem 3.1.2]{HRW25}. Moreover, $(f_0)_*$ is surjective with kernel generated by 
	\[v_i-c_iu^{p^i-1},\qquad y_{\delta}-c_{\delta}u^{|y_{\delta}|/2}\]
	for $i\geq1$, where $c_i=u_i$ for $1\leq i\leq m-1$.
\end{lem}
\begin{proof}
	By \cite[Remark 2.5]{BSY22} $\mathbb{W}(L)$ is flat over $\Ss_p$. Hence there is an isomorphism
	\begin{equation*}
		\begin{split}
			& \pi_*(\MU\otimes\mathbb{W}(L)\otimes\Sigma^{\infty}_+(W_0^{\times m-1}\times W_2)/p^k) \\
			& \cong\pi_*(\mathbb{W}(L)\otimes_{\Ss_p}(\Ss_p\otimes\MU\otimes\Sigma^{\infty}_+(W_0^{\times m-1}\times W_2)/p^k)) \\
			& \cong W(L)/p^k[x_1^{\pm},\cdots,x_{m-1}^{\pm},u,v_i,y_{\delta}],
		\end{split}
	\end{equation*}
	for $k\geq1$, by \cite[Theorem 3.1.2]{HRW25}. In particular, $\pi_*(\MU\otimes\mathbb{W}(L)\otimes\Sigma^{\infty}_+(W_0^{\times m-1}\times W_2)/p^k)$ concentrates in even degrees. Therefore, we obtain that
	\begin{equation*}
		\begin{split}
			& \pi_*(\MU\otimes\mathbb{W}(L)\otimes\Sigma^{\infty}_+(W_0^{\times m-1}\times W_2)/(p^k,(x_1-1)^{i_1},\cdots,(x_{m-1}-1)^{i_{m-1}}))\\
			& \cong W(L)/p^k[u_1,\cdots,u_{m-1},u,v_i,y_{\delta}]/(u_1^{i_1},\cdots,u_{m-1}^{i_{m-1}}),
		\end{split}
	\end{equation*}
	for $i_1,\cdots,i_{m-1}>0$. In particular, the $\mathrm{lim}^1$-term in the Milnor exact sequence for computing $\pi_*(W(0,m))$ vanishes, which gives the desired isomorphism. Since the positive degree of $\pi_*(e_m(L))$ is spanned by $cu^k$ for $k\geq 1$ and $c\in \pi_0(e_m(L))$, the kernel also admits desired generators by our choices.
\end{proof}

Note that there are only finitely many $y_{\delta}$ in each degree. Since $v_i-c_iu^{p^i-1}\in\mathrm{Ker}((f_0)_*)$, the suspension map in \cite[Example A.2.4]{HW22} gives the element 
\[z_i:=\sigma^2(v_i-c_iu^{p^i-1})\in \pi_*(T\otimes_{W(0,m)}e_m(L))\]
for $i\geq1$. The following lemma is similar to \cite[Lemma 2.3]{Zha26}, which will be used to show that $K(m+n)\otimes \Wn^{hT^n}\neq0$.

\begin{lem}\label{lem:zi}
	Suppose that $m\geq2$, or that $m=1$ and $p$ is odd. In the quotient ring
	\[\pi_*(T\otimes_{W(0,m)}e_m(L))/(p,u_1,\cdots,u_{m-1},u),\]
	we have that $\overline{z}_i\doteq (\overline{z}_1)^{p^{i-1}}$, where $\overline{z}_i$ is the image of $z_i$. In particular, it is a free $L[\overline{z}_1]$-module.
\end{lem}
\begin{proof}
	By \cite[Proposition 3.6]{Ang08}, Lemma \ref{lem:pi(W0)} implies that there is an isomorphism 
	\[\pi_*(e_m(L)\otimes_{W(0,m)}e_m(L))\cong\Lambda_{\pi_*(e_m(L))}(\sigma(v_i-c_iu^{p^i-1}),\sigma(y_{\delta}-c_{\delta}u^{|y_{\delta}|/2})).\]
	Using the Tor spectral sequence, we obtain that there is a filtration on $\pi_*(T\otimes_{W(0,m)}e_m(L))$ such that
	\begin{equation}\label{eq:gr/p}
		\mathrm{gr}(\pi_*(T\otimes_{W(0,m)}e_m(L)))\cong\Gamma_{\pi_*(e_m(L))}(z_i,\sigma^2(y_{\delta}-c_{\delta}u^{|y_{\delta}|/2})).
	\end{equation}
	In particular, $\pi_*(T\otimes_{W(0,m)}e_m(L))$ is even and $(p,u_1,\cdots,u_{m-1},u)$ is a regular sequence. Hence there is an isomorphism $\pi_*((T\otimes_{W(0,m)}e_m(L))\otimes_{e_m(L)}L)\cong\pi_*(T\otimes_{W(0,m)}e_m(L))/(p,u_1,\cdots,u_{m-1},u)$.
	
	Now, applying $T\otimes_{W(0,m)}-$ to the canonical $\Ee_{\infty}$-map $e_m(L)\to L$ and base change along it gives an $\Ee_{\infty}$-map
	\[f\colon(T\otimes_{W(0,m)}e_m(L))\otimes_{e_m(L)}L\longrightarrow T\otimes_{W(0,m)}L.\]
	Using the isomorphism 
	\[\pi_*(W(0,m))\cong W(L)[[u_1,\cdots,u_{m-1}]][u,v_i,y_{\delta}]\cong W(L)[[u_1,\cdots,u_{m-1}]][u,v_i-c_iu^{p^i-1},y_{\delta}-c_{\delta}u^{|y_{\delta}|/2}],\]
	the same computation implies that there is a filtration on $\pi_*(T\otimes_{W(0,m)}L)$ such that
	\[\mathrm{gr}(\pi_*(T\otimes_{W(0,m)}L))\cong\Gamma_L(\sigma^2p,\sigma^2(u_1),\cdots,\sigma^2(u_{m-1}),\sigma^2(u),z_i,\sigma^2(y_{\delta}-c_{\delta}u^{|y_{\delta}|/2})).\]
	Moreover, $\mathrm{gr}(f_*)$ sends $\overline{z}_i$ and $\overline{\sigma^2(y_{\delta}-c_{\delta}u^{|y_{\delta}|/2})}$ in $\mathrm{gr}(\pi_*(T\otimes_{W(0,m)}e_m(L))/(p,u_1,\cdots,u_{m-1},u))$ to $z_i$ and $\sigma^2(y_{\delta}-c_{\delta}u^{|y_{\delta}|/2})$ in $\mathrm{gr}(\pi_*(T\otimes_{W(0,m)}L))$, respectively, where $\overline{\sigma^2(y_{\delta}-c_{\delta}u^{|y_{\delta}|/2})}$ is the image of $\sigma^2(y_{\delta}-c_{\delta}u^{|y_{\delta}|/2})$. Since both associated gradeds are finite in each degree, $f_*$ is injective in homotopy groups. Hence it suffices to show that $z_i\doteq z_1^{p^{i-1}}$ in $\pi_*(T\otimes_{W(0,m)}L)$. When $p^i-1>1$, we have that 
	\[z_i=\sigma^2(v_i-c_iu^{p^i-1})=\sigma^2(v_i)-c_iu^{p^i-2}\cdot\sigma^2(u)=\sigma^2(v_i)\]
	for $i\geq1$, since $\sigma^2$ is $\pi_*(W(0,m))$-linear and $u$ is sent to zero in $L$. When $p^i-1=1$, $c_1=u_1$ is killed in $L$, since $m\geq2$. Hence we also have $z_i=\sigma^2(v_i)$. Moreover, for each $i\geq1$, $\sigma^2(v_i)$ lies in the image of $\pi_*(T\otimes_{\MU}\mathbb{F}_p)\to\pi_*(T\otimes_{W(0,m)}L)$. The desired formula now follows from \cite[Proposition 2.4.2]{HW22} and the argument in \cite[Lemma 2.3]{Zha26}. In particular, $(\overline{z}_1)^k$ is nonzero for all $k\geq0$, which implies that $\pi_*(T\otimes_{W(0,m)}e_m(L))/(p,u_1,\cdots,u_{m-1},u)$ is a free $L[\overline{z}_1]$-module by the isomorphism \eqref{eq:gr/p}.
\end{proof}

When $m=1$ and $p=2$, $|v_1|=|u|=2$. Hence, the $\Ee_{\infty}$-map
\[f^{\prime}_0\colon(\mathbb{W}(L)\otimes\MU)_2\longrightarrow e_1(L)\]
is already surjective on homotopy groups. Without loss of generality, we may assume that $v_1$ is sent to $u$. Consequently, a similar computation as in Lemma \ref{lem:pi(W0)} implies that $v_i-c_iv_1^{2^i-1}$ belongs to the kernel of $(f^{\prime}_0)_*$ for some $c_i\in W(L)$ and $i\geq2$. The element $\sigma^2(v_i-c_iv_1^{2^i-1})\in\pi_*(T\otimes_{(\mathbb{W}(L)\otimes\MU)_2}e_1(L))$ is also denoted by $z_i$ for $i\geq2$. The same argument as the proof of Lemma \ref{lem:zi} gives the following lemma.

\begin{lem}\label{lem:wi}
	Suppose that $m=1$ and $p=2$. In the quotient ring
	\[\pi_*(T\otimes_{(\mathbb{W}(L)\otimes\MU)_2}e_1(L))/(2,u),\]
	we have that $\overline{z}_i\doteq (\overline{z}_2)^{2^{i-2}}$, where $\overline{z}_i$ is the image of $z_i$. In particular, it is a free $L[\overline{z}_2]$-module.
\end{lem}

Now, we define $\Wn$ inductively and establish some of its properties. Given a sequence $\{x_{\alpha}\}$ in a graded ring $A$, we call $\{x_{\alpha}\}$ finitely regular if it satisfies the following two properties 
\begin{enumerate}
	\item[$(1)$] each $x_{\alpha}$ has positive even degree, and there are only finitely many $x_{\alpha}$ in any positive even degree,
	\item[$(2)$] every finite subfamily of $\{x_{\alpha}\}$ is a regular sequence;
\end{enumerate}
compare with \cite[Lemma 3.1]{Zha26}. We begin with the following proposition.

\begin{prop}\label{prop:ind}
	Let $f\colon W\to R$ be a map of $p$-local connective even $\Ee_{\infty}$-rings such that $f_*$ is surjective. Suppose that 
	\begin{enumerate}
		\item[$(1)$] its kernel is generated by a finitely regular sequence $\{h_{\delta}\}$, 
		\item[$(2)$] $\pi_*(R)$ is isomorphic to $\pi_0(W)[x_{\alpha}]/(f_{\beta})$ with $|x_{\alpha}|>0$ such that there are only finitely many $x_{\alpha}$  in any positive even degree, and $\{f_{\beta}\}$ is a finitely regular sequence in $\pi_0(W)[x_{\alpha}]$.
	\end{enumerate}
	Then there is an isomorphism
	\[\pi_*(T\otimes_WR)\cong\frac{\pi_0(W)[x_{\alpha},X_{\delta,i}\ |\ i\geq0]}{(f_{\beta},F_{\delta,i})},\]
	where $X_{\delta,i}$ is a lift of $\gamma_{p^i}(\sigma^2(h_{\delta}))$ in the spectral sequence 
	\[E^2\cong\mathrm{Tor}^{\Lambda_{\pi_*(R)}(\sigma(h_{\delta}))}(\pi_*(R),\pi_*(R))\cong\Gamma_{\pi_*(R)}(\sigma^2(h_{\delta}))\Longrightarrow\pi_*(T\otimes_WR),\]
	and $F_{\delta,i}=X_{\delta,i}^{p}-p\epsilon_iX_{\delta,i+1}-G_{\delta,i}$ with $\epsilon_i=\frac{(p^{i+1})!}{p((p^i)!)^p}\in\mathbb Z_{(p)}^\times$ for some $G_{\delta,i}\in\pi_0(W)[x_{\alpha},X_{\delta,i}]$. In particular, for each $k>0$, there are only finitely many $X_{\delta,i}$ of degree $2k$. Moreover, $\{f_{\beta},F_{\delta,i}\}$ is a finitely regular sequence in $\pi_0(W)[x_{\alpha},X_{\delta,i}]$. 
\end{prop}
\begin{proof}
	The proof is the same as that of \cite[Proposition 3.2]{Zha26}.
\end{proof}

We take 
\[W:=\begin{cases}
	W(0,m), & \text{$m\geq2$ or $m=1$ and $p$ is odd} \\
	(\mathbb{W}(L)\otimes\MU)_2, & \text{$m=1$ and $p=2$}.
\end{cases}\]
in the result of this paper. For $1\leq i\leq n$, the inclusion of $i$-th $T$ into $T^n$ gives a map $\iota_i\colon T\otimes_{W}e_m(L)\to T^n\otimes_{W}e_m(L)$. We denote the image of $z_j\in\pi_*(T\otimes_{W}e_m(L))$ under $(\iota_i)_*$ by $\zeta_{i,j}$ for $j\geq1$. For convenience, we also write $\zeta_i$ for $\zeta_{i,1}$ when $m\geq2$ or $m=1$ and $p$ is odd, and for $\zeta_{i,2}$ when $m=1$ and $p=2$. Now, we can construct the $\Ee_{\infty}$-$W$-algebra $\Wn$ by induction on $n$. The construction of $\Wn$ is similar to \cite[Proposition 3.4]{Zha26}.

\begin{prop}\label{prop:wn}
	For $n\geq1$, there is an even $\Ee_{\infty}$-$W$-algebra $\Wn$ equipped with a $T^n$-action and a $T^n$-equivariant $\Ee_{\infty}$-$W$-algebra map $g_n\colon\TT{n}_We_m(L)\to \Wn$. Moreover, $\pi_*(\Wn)$ is a free $\pi_*(e_m(L))$-module and $\pi_*(\Wn)/(p,u_1,\cdots,u_{m-1},u)$ is a free $L[\zeta_1,\cdots,\zeta_n]$-module. Here, for each $1\leq i\leq n$, we consider $\zeta_i$ as an element of $\pi_*(\Wn)/(p,u_1,\cdots,u_{m-1},u)$ via $(g_n)_*$.
\end{prop}
\begin{proof}
	First, define $W(1,m):=T\otimes_W e_m(L)$ and take $g_1=id$. For a fixed $n\geq1$, suppose that we have constructed $\Wn$ with the desired properties, and 
	\[\pi_*(\Wn)\cong\pi_*(e_m(L))[x_{\alpha}]/(f_{\beta}(x_{\alpha})),\]
	where there are only finitely many $x_{\alpha}$ in any positive even degree and $\{f_{\beta}\}$ is finitely regular. Note that $W(1,m)$ satisfies the inductive hypothesis by Lemma \ref{lem:pi(W0)} and Proposition \ref{prop:ind}. 
	
	Since $\Wn$ is even, the homotopy fixed point spectral sequence for $\pi_*(\Wn^{hT^n})$ collapses at the $E_2$-page. Hence $\Wn^{hT^n}$ is even and the canonical map $\pi_*(\Wn^{hT^n})\to\pi_*(\Wn)$ is surjective. Consequently, we can choose a lift of $x_{\alpha}\in\pi_*(\Wn)$ in $\pi_*(\Wn^{hT^n})$, which is still denoted by $x_{\alpha}$. Since $|x_{\alpha}|$ is positive even and $W$ is a connective even $\Ee_{\infty}$-ring, there is an $\Ee_{\infty}$-$W$-map $W\otimes W_{|x_{\alpha}|}\to \Wn^{hT^n}$ which sends the generator $a_{\alpha}\in\pi_{|x_{\alpha}|}(W\otimes W_{|x_{\alpha}|})$ to $x_{\alpha}$ by \cite[Proposition 3.1.5]{HRW25}. Set $X=\prod_{\alpha}W_{|x_{\alpha}|}$. Since there are only finitely many $x_{\alpha}$ in any positive even degree, there is an isomorphism $W_*(X)\cong \pi_*(W)[a_{\alpha},b_{\alpha,\gamma(\alpha)}]$ with $|b_{\alpha,\gamma(\alpha)}|>|a_{\alpha}|$, by \cite[Theorem 3.1.2]{HRW25}. Moreover, multiplying above maps together and using the adjunction between the fixed point functor and the forgetful functor, we obtain a $T^n$-equivariant $\Ee_{\infty}$-$W$-algebra map
	\[f_n\colon W\otimes X\simeq \mathrm{colim}_kW\otimes \prod_{|x_{\alpha}|<2k}W_{|x_{\alpha}|}\longrightarrow \Wn,\]
	which sends $a_{\alpha}\in W_*(X)$ to $x_{\alpha}\in\pi_*(\Wn)$. 
	
	We complete the induction for $m\geq2$ or $m=1$ and $p$ is odd. The case $m=1$ and $p=2$ is similar. In this case, we have that
	\[W_*(X)\cong\pi_*(e_m(L))[v_i,y_{\delta},a_{\alpha},b_{\alpha,\gamma(\alpha)}]\cong\pi_*(e_m(L))[v_i-c_iu^{p^i-1},y_{\delta}-c_{\delta}u^{|y_{\delta}|/2},a_{\alpha},b_{\alpha,\gamma(\alpha)}].\]
	Since $f_n$ is $W$-linear, the kernel of $(f_n)_*$ is generated by 
	\[\{v_i-c_iu^{p^i-1},y_{\delta}-c_{\delta}u^{|y_{\delta}|/2},f_{\beta}(a_{\alpha}),F_{\alpha,\gamma}(a_{\alpha})-b_{\alpha,\gamma(\alpha)}\},\]
	where $F_{\alpha,\gamma}(x_{\alpha})$ is the image of $b_{\alpha,\gamma(\alpha)}$. In particular, every finite subfamily is a regular sequence. Define $W(n+1,m):=T\otimes_{W\otimes X}\Wn$ and 
	\[g_{n+1}\colon\TT{n+1}_We_m(L)\simeq T\otimes_W(\TT{n}_We_m(L))\xrightarrow{T\otimes_Wg_n}T\otimes_W\Wn\xrightarrow{T\otimes_{\eta}\Wn}W(n+1,m),\]
	where $\eta\colon W\to W\otimes X$ is the unit map. Considering $v_i-c_iu^{p^i-1}$ as an element of $\mathrm{Ker}(f_n)_*$ for $i\geq1$, the $W$-linearity implies that the element $\sigma^2(v_i-c_iu^{p^i-1})\in\pi_*(W(n+1,m))/(p,u_1,\cdots,u_{m-1},u)$ given by the suspension map is exactly $\zeta_{n+1,i}$ defined above. Using Lemma \ref{lem:zi}, we obtain the remaining properties by an argument analogous to that in the proof of \cite[Proposition 3.4]{Zha26}.
\end{proof}

\section{The height of $\Wn^{hT^n}$ and the proof of Theorem \ref{thm1}}\label{sec4}
In this section, we show that the height of $\Wn^{hT^n}$ is $m+n$. This allows us to determine the height of $(\TT{n}e_m(L))^{hT^n}$, which is a key step in the proof of Theorem \ref{thm1}. 

For this, note that, using the functor $\mathrm{sh}^k$ in \cite[Variant 3.2]{ABM26}, Construction 3.3 and Proposition 3.12 in \cite{ABM26} extend to $k<0$. Hence, for integers $w$ and $k$, there is an $\Ee_2$-ring $\Ss[x_{2k,w}]$ in $\mathrm{Fun}(\Zz^{\mathrm{ds}},\mathrm{Sp})$ whose underlying $\Ee_1$-ring is the free $\Ee_1$-ring generated by $\Sigma^{2k}\Ss$ concentrated in $w$. Here $\mathrm{Fun}(\Zz^{\mathrm{ds}},\mathrm{Sp})$ is considered as a symmetric monoidal category via the Day convolution. We call $w$ the weight of $x_{2k,w}$, denoted by $w(x_{2k,w})=w$. Moreover, there is an $\Ee_2$-map $\mathrm{Free}_{\Ee_2}(\Sigma^{2k}\Ss(w))\to \Ss[x_{2k,w}]$ which admits a cell decomposition with cells concentrated in even degrees, where $\Ss(w)$ denotes the associated spectrum with $\Ss$ placed in weight $w$. We also need the following lemma.

\begin{lem}\label{lem:complete}
	Let $M$ be a complete filtered spectrum. If $\pi_*(\mathrm{gr}(M))$ concentrates in even degrees and is nonzero, then $\mathrm{colim}(M)$ is nonzero. 
\end{lem}
\begin{proof}
	Since $M$ is complete, we have that
	\[M_s\simeq\lim_t M_s/M_{s+t}.\]
	For $t>0$, using the fiber sequence
	\[\mathrm{gr}_{s+t}(M)\to M_s/M_{s+t+1}\to M_s/M_{s+t},\]
	we see that $\pi_*(M_s/M_{s+t})$ concentrates in even degrees and $\pi_*(M_s/M_{s+t+1})\to\pi_*(M_s/M_{s+t})$ is surjective. Hence the $\mathrm{lim}^1$ term in the Milnor exact sequence for computing $\pi_*(M_s)$ vanishes, and $\pi_*(M_s)$ is even. This implies that $\pi_*(M_s)\to\pi_*(M_{s-1})$ is injective and $\pi_*(M_s)\to\pi_*(\mathrm{gr}_s(M))$ is surjective. Therefore, $\mathrm{colim}(M)$ is nonzero.
\end{proof}

\begin{prop}\label{prop:wn nonzero}
	For $n\geq1$, $K(m+n)\otimes\Wn^{hT^n}\not\simeq0$.
\end{prop}
\begin{proof}
	First, suppose that $m\geq2$ or $m=1$ and $p$ is odd. Let $T_i$ be the $i$-th $T$ in $T^n$. Since $\Wn$ is connective even, multiplying by $t_i\in\pi_{-2}(\Wn^{hT_i})$ and taking $T^n/T_i$-fixed point induce a cofiber sequence
	\[\Sigma^{-2}W(n,m)^{hT^n}\xrightarrow{\ t_i\ }W(n,m)^{hT^n}\longrightarrow W(n,m)^{hT^n/T_i}\]
	of even spectra, for $1\leq i\leq n$. Iterating this $n$ times and using Proposition \ref{prop:wn}, we obtain that
	\begin{equation}\label{eq:/t1tn}
		\pi_*(W(n,m)^{hT^n})/(t_1,\cdots,t_n)\cong\pi_*(W(n,m))
	\end{equation}
	is a free $L[\zeta_1,\cdots,\zeta_n]$-module after modulo $(p,u_1,\cdots,u_{m-1},u)$. By \cite[Lemma A.4.1]{HW22} and the argument in \cite[Proposition 4.3]{Zha26}, we see that
	\[v_j-c_ju^{p^j-1}\equiv\sum_{i=1}^nt_i\zeta_{i,j}\mod (t_1,\cdots,t_n)^2,\]
	where we consider $v_j-c_ju^{p^j-1}$ as an element in $\pi_*(\Wn^{hT^n})$ via the $W$-algebra structure on $\Wn$ for $j\geq1$. By the proof of Proposition \ref{prop:wn} and the isomorphism \eqref{eq:/t1tn}, $v_j-c_ju^{p^j-1}$ belongs to the kernel of the quotient map
	\[\pi_*(W(n,m)^{hT^n})\longrightarrow \pi_*(W(n,m)^{hT^n})/(t_1,\cdots,t_n)\cong\pi_*(W(n,m)).\]
	Consequently, we obtain that
	\[v_j-c_ju^{p^j-1}=\sum_{i=1}^nt_ib_{i,j}\in\pi_*(W(n,m)^{hT^n}),\]
	for some $b_{i,j}\in\pi_{2p^j}(W(n,m)^{hT^n})$. By \cite[Lemma 4.2]{Zha26}, $\{t_1,\cdots,t_n\}$ is a regular sequence in $\pi_*(W(n,m)^{hT^n})/(p,u_1,\cdots,u_{m-1},u)$. Hence Lemma \ref{lem:zi} implies that $\overline{b}_{i,j}\doteq\zeta_i^{p^{j-1}}$ in 
	\begin{equation}\label{eq:A}
		A:=\pi_*(\Wn^{hT^n})/(p,u_1,\cdots,u_{m-1},u,t_1,\cdots,t_n)\cong\pi_*(W(n,m))/(p,u_1,\cdots,u_{m-1},u)
	\end{equation}
	by the isomorphism \eqref{eq:/t1tn}, where $\overline{b}_{i,j}$ is the image of $b_{i,j}$ in $A$. 
	
	We show that this implies that $K(m+n)\otimes\Wn^{hT^n}$ is nonzero. The cell decomposition mentioned in the beginning of this section implies that there is an $\Ee_2$-map 
	\[\Ss[U,U_1,\cdots,U_{m-1},T_1,\cdots,T_n]\to\Wn^{hT^n},\]
	which sends $U$ to $u$, $U_i$ to $u_i$, and $T_j$ to $t_j$. Base change along the map $\Ss\to\Wn^{hT^n}$, we obtain an $\Ee_2$-$\Wn^{hT^n}$-algebra map $B:=\Wn^{hT^n}[U,U_1,\cdots,U_{m-1},T_1,\cdots,T_n]\to\Wn^{hT^n}$. This gives a $B$-action on $\Wn^{hT^n}$. We write $\Wn^{hT^n}_B$ for $\Wn^{hT^n}$ equipped with this $B$-action. Consider $B$ as an object of $\mathrm{Fun}(\Zz^{\mathrm{ds}},\mathrm{Mod}_{\Wn^{hT^n}})$ by defining
	\[w(U)=1,\qquad w(U_i)=p^m-p^i,\qquad w(T_j)=p^m-1.\] 
	Now, left Kan extension along the symmetric monoidal functor $i\colon \Zz^{\mathrm{ds}}\to\Zz_{>}$ induces a symmetric monoidal functor $i_!\colon\mathrm{Fun}(\Zz^{\mathrm{ds}},\mathrm{Mod}_{\Wn^{hT^n}})\to\mathrm{Fun}(\Zz_{>},\mathrm{Mod}_{\Wn^{hT^n}})$. Set $C:=i_!(B)$. We have that $C_0=B$ as $\Wn^{hT^n}$-modules. Let $i_0\colon\{*\}\to\Zz_{>}$ be the inclusion of zero. Since $i_0$ is also symmetric monoidal, we have that
	\[\mathrm{Map}_{\mathrm{CAlg}(\mathrm{Fun}(\Zz_{>},\mathrm{Mod}_{\Wn^{hT^n}}))}((i_0)_!B,C)\simeq \mathrm{Map}_{\mathrm{CAlg}(\mathrm{Mod}_{\Wn^{hT^n}})}(B,B),\]
	where $(i_0)_!$ is the left Kan extension along $i_0$. Here $B$ admits an $\Ee_{\infty}$-ring structure by \cite[Remark 4.5]{RSS26}. In particular, there is an $\Ee_{\infty}$-map $(i_0)_!B\to C$ corresponding to the identity map of $B$. Since 
	\[\mathrm{gr}_0(C)\simeq\mathrm{cofib}(C_1\to C_0)=\mathrm{cofib}(\bigoplus_{w\geq1}B_w\to B)\simeq \Wn^{hT^n},\]
	the $B\simeq\mathrm{gr}((i_0)_!B)$-action on $\mathrm{gr}(C)$ is induced by the augmentation $B\to \Wn^{hT^n}$. Hence we obtain that
	\begin{equation*}
		\begin{split}
			& \mathrm{gr}((i_0)_!(\Wn^{hT^n}_B)\otimes_{(i_0)_!B}C)\simeq\Wn^{hT^n}\otimes_{B}\mathrm{gr}(C) \\
			& \simeq(\Wn^{hT^n}_B\otimes_{B}\Wn^{hT^n})\otimes_{\Wn^{hT^n}}\mathrm{gr}(C),
		\end{split}
	\end{equation*}
	where $B$ acts on the second $\Wn^{hT^n}$ by the augmentation. Since $\mathrm{gr}(C)\simeq B$, we see that 
	\begin{equation}\label{eq:grC}
		\pi_*(\mathrm{gr}((i_0)_!(\Wn^{hT^n}_B)\otimes_{(i_0)_!B}C)/p)\cong A[U,U_1,\cdots,U_{m-1},T_1,\cdots,T_n]
	\end{equation}
	by the isomorphism \eqref{eq:A}. Note that, since the map $B\to\Wn^{hT^n}$ is an $\Ee_2$-$\Wn^{hT^n}$-algebra map, $(i_0)_!(\Wn^{hT^n}_B)\otimes_{(i_0)_!B}C)$ is a filtered $\Ee_1$-$\Wn^{hT^n}$-algebra and $(i_0)_!(\Wn^{hT^n}_B)\otimes_{(i_0)_!B}C)/p$ admits a $\Wn^{hT^n}$-module structure.
	
	We construct a $\Wn^{hT^n}$-module $N$ from $(i_0)_!(\Wn^{hT^n}_B)\otimes_{(i_0)_!B}C/p$ such that $K(m+n)\otimes N$ is nonzero. Since $N$ is a retract of $\Wn^{hT^n}\otimes N$, this shows that $K(m+n)\otimes\Wn^{hT^n}$ is also nonzero. For this, consider the following elements in $\pi_*(B)$
	\[V_j=\begin{cases}
		U_jU^{p^j-1}+\displaystyle\sum_{i=1}^nT_ib_{i,j}, & 1\leq j<m,\\
		c_jU^{p^j-1}+\displaystyle\sum_{i=1}^nT_ib_{i,j}, & m\leq j\leq m+n.
	\end{cases}\]
	Then the equation $\overline{b}_{i,j}\doteq\zeta_i^{p^{j-1}}$ in $A$ implies that the image of the weight $p^m-1$ summand of $V_j$ in $A[U,U_1,\cdots,U_{m-1},T_1,\cdots,T_n]$ is
	\[V^{\prime}_j=\begin{cases}
		U_jU^{p^j-1}+a_j\displaystyle\sum_{i=1}^nT_i\zeta_i^{p^{j-1}}, & 1\leq j<m,\\
		c_mU^{p^m-1}+a_m\displaystyle\sum_{i=1}^nT_i\zeta_i^{p^{m-1}}, & j=m,\\
		a_j\displaystyle\sum_{i=1}^nT_i\zeta_i^{p^{j-1}}, & m<j\leq m+n,
	\end{cases}\]
	where $a_j$ are units in $L$. Let $D=\mathrm{det}(\zeta_i^{p^{j-1}})_{1\leq i,j\leq n}\in A.$ Since $A$ is a free $L[\zeta_1,\cdots,\zeta_n]$-module by Proposition \ref{prop:wn}, $A[D^{-1}]$ is nonzero. We claim that $(\overline{V}_1^{\prime},\cdots,\overline{V}_{m+n}^{\prime})$ is a regular sequence in 
	\[\pi_*(B)/(p,u_1,\cdots,u_{m-1},u,t_1,\cdots,t_n)[D^{-1},U^{-1}]\cong A[D^{-1},U^{\pm},U_1,\cdots,U_{m-1},T_1,\cdots,T_n].\]
	Indeed, since $U$ is invertible in this ring, it suffices to show that $(\overline{V}_m^{\prime},\cdots,\overline{V}_{m+n}^{\prime})$ is a regular sequence in $A[D^{-1},U^{\pm},T_1,\cdots,T_n]$. Moreover, we have that
	\[\begin{pmatrix}
		\overline{V}_{m}^{\prime}-c_mU^{p^m-1}\\
		\vdots\\
		\overline{V}_{m+n-1}^{\prime}
	\end{pmatrix}
	=
	\begin{pmatrix}
		a_m\zeta_1^{p^{m-1}} & \cdots & a_m\zeta_n^{p^{m-1}} \\
		\vdots & \ddots & \vdots \\
		a_{m+n-1}\zeta_1^{p^{m+n-2}} & \cdots & a_{m+n-1}\zeta_n^{p^{m+n-2}}
	\end{pmatrix}
	\begin{pmatrix}
		T_1\\
		\vdots\\
		T_n
	\end{pmatrix}.\]
	Since $p=0$ in $A$, we have that
	\[\mathrm{det}\begin{pmatrix}
		a_m\zeta_1^{p^{m-1}} & \cdots & a_m\zeta_n^{p^{m-1}} \\
		\vdots & \ddots & \vdots \\
		a_{m+n-1}\zeta_1^{p^{m+n-2}} & \cdots & a_{m+n-1}\zeta_n^{p^{m+n-2}}
	\end{pmatrix}
	=\left(\prod_{i=m}^{m+n-1}a_i\right)D^{p^{m-1}}\]
	is invertible in $A[D^{-1},U^{\pm},T_1,\cdots,T_n]$. Hence we see that
	\[A[D^{-1},U^{\pm},T_1,\cdots,T_n]/(\overline{V}_m^{\prime},\cdots,\overline{V}_{m+n-1}^{\prime})\cong A[D^{-1},U^{\pm}].\]
	It remains to show that the image of $\overline{V}_{m+n}^{\prime}$ in this quotient is a nonzerodivisor. Let
	\[\zeta^i=(\zeta_1^{p^{m-1+i}},\cdots,\zeta_n^{p^{m-1+i}}),\qquad T=(T_1,\cdots,T_n)^T\]
	for $0\leq i\leq n$. Then $\overline{V}_m^{\prime}=a_m\zeta^0\cdot T+c_mU^{p^m-1}$ and $\overline{V}_{m+i}^{\prime}=a_{m+i}\zeta^i\cdot T$ in $A[D^{-1},U^{\pm},T_1,\cdots,T_n]$ for $1\leq i\leq n$, which implies that
	\begin{equation}
		a_m\zeta^0\cdot T+c_mU^{p^m-1}=a_{m+i}\zeta^i\cdot T=0
	\end{equation}
	in $A[D^{-1},U^{\pm}]$, for $1\leq i\leq n-1$. We have that
	\[\mathrm{det}\begin{pmatrix}
		\zeta^0\\
		\vdots\\
		\zeta^{n-1}
	\end{pmatrix}
	=D^{p^{m-1}},\]
	is a unit in $A[D^{-1},U^{\pm}]$, as $p=0$. In particular, there is a unique sequence $d_0,\cdots,d_{n-1}$ in $A[D^{-1},U^{\pm}]$ such that $\zeta^n=d_0\zeta^0+\cdots+d_{n-1}\zeta^{n-1}$, which implies that $\overline{V}_{m+n}^{\prime}=-a_m^{-1}a_{m+n}c_md_0U^{p^m-1}$. Now, we have that
	\[D^{p^m}=\mathrm{det}\begin{pmatrix}
		\zeta^1\\
		\vdots\\
		\zeta^n
	\end{pmatrix}
	=\sum_{i=0}^{n-1}d_i\cdot\mathrm{det}\begin{pmatrix}
		\zeta^1\\
		\vdots\\
		\zeta^i
	\end{pmatrix}
	=d_0\cdot\mathrm{det}\begin{pmatrix}
		\zeta^1\\
		\vdots\\
		\zeta^0
	\end{pmatrix}
	=(-1)^{n-1}d_0D^{p^{m-1}}.\]
	Hence $d_0=(-1)^{n-1}D^{p^m-p^{m-1}}$ is a unit in $A[D^{-1},U^{\pm}]$, and so is $\overline{V}_{m+n}^{\prime}$.
	
	Choose a lift of $D$ in $\pi_*(\Wn^{hT^n})$, which is denoted by $\widetilde{D}$. Let 
	\[M:=\left((i_0)_!(\Wn^{hT^n}_B)\otimes_{(i_0)_!B}C/p\right)[U^{-1},\widetilde{D}^{-1},V_{m+n}^{-1}]\]
	be a filtered $\Wn^{hT^n}$-module and $\widehat{M}$ be its completion, given by $\widehat{M}_s=\mathrm{cofib}(\lim_t M_t\to M_s)$. Set
	\[\mathcal{N}:=\widehat{M}/(V_1,\cdots,V_{m+n-1}).\]
	By the equivalence \eqref{eq:grC} and the discussion in the preceding paragraph, we have that $\pi_*(\mathrm{gr}(\mathcal{N}))\cong A[U^{\pm},\widetilde{D}^{-1}]$, which concentrates in even degrees. By construction, $\mathcal{N}$ is complete. Lemma \ref{lem:complete} then implies that $N:=\mathrm{colim}\mathcal{N}_s$ is nonzero. 
	
	Finally, by our choice, the actions of $V_j$ and $v_j$ on $\pi_*(\Wn^{hT^n}_B)$ are the same. Hence, for $1\leq j\leq m+n$, $X_j:=V_j-v_j$ acts on $\pi_*(\Wn^{hT^n}_B)$ trivially. Now, as a $\pi_*(B)$-module, the projective dimension of $\pi_*(\Wn^{hT^n}_B)$ is $m+n$, which implies that $X_j^{m+n+1}$ acts on $\Wn^{hT^n}_B$ trivially in the category of $B$-modules by \cite[Proposition 1.4.1]{HL13}. Hence it also acts on $N$ trivially. Moreover, it is well-known that, for any $\Ee_1$-ring $R$ and $x\in\pi_*(R)$, $x^2$ acts trivially on $R/x$. Hence, for $1\leq j\leq m+n-1$, we obtain that the actions of $p^2$ and $V_j^2$ on $N$ are trivial, and
	\begin{equation*}
		\begin{split}
			v_j^{m+n+2}=(V_j-X_j)^{m+n+2}= & (-1)^{m+n+2}X_j^{m+n+2}+(-1)^{m+n+1}(m+n+2)V_jX_j^{m+n+1}\\
			& +\sum_{i=2}^{m+n+2}(-1)^{m+n+2-i}\binom{m+n+2}{i} V_j^iX_j^{m+n+2-i}=0.
		\end{split}
	\end{equation*}
	In particular, $N/(p,v_1,\cdots,v_{m+n-1})$ is nonzero. On the other hand, since $V_{m+n}$ acts invertibly on $N$, $v_{m+n}=V_{m+n}-X_{m+n}$ is invertible in $N$. Hence $N/(p,v_1,\cdots,v_{m+n-1})[v_{m+n}^{-1}]$ is nonzero. By \cite[Lemma 4.4]{Zha26}, $K(m+n)\otimes N$ is nonzero. 
	
	The proof of the case when $m=1$ and $p=2$ is the same as that of \cite[Theorem 1.3]{Zha26} with $\ell$ replaced by $e_1(L)$.
\end{proof}

\begin{proof}[Proof of Theorem \ref{thm1}]
	First, we assume that $m\geq1$. The $\Ee_{\infty}$-map \eqref{eq:forget} and Proposition \ref{prop:wn} imply that there is an $\Ee_1$-map 
	\[K(m+n)\otimes(\TT{n}e_m(L))^{hT^n}\to K(m+n)\otimes\Wn^{hT^n}.\]
	Since $K(m+n)\otimes\Wn^{hT^n}$ is nonzero by Proposition \ref{prop:wn nonzero}, so is $K(m+n)\otimes(\TT{n}e_m(L))^{hT^n}$. By the argument before Lemma \ref{lem:e3}, this finishes the proof of Theorem \ref{thm1} when $m\geq1$. 
	
	Next, we consider the case $m=-1$. Since $R$ is connective, $K^{(n)}(R)$ has height $n-1$ by Theorem E and Remark 9.13 in \cite{BSY22}. Hence the same argument as in the introduction implies that the height of $(T^n\otimes R)^{hT^n}$ is at most $n-1$. It remains to show that $K(n-1)\otimes(T^n\otimes R)^{hT^n}$ is nonzero. Since $R$ has height $-1$, $\mathbb{Q}\otimes R\simeq R[1/p]\simeq0$, which implies that $p^N$ is zero in $\pi_0(R)$ for some $N>0$. In particular, there is a map of commutative rings $\pi_0(R)\to(\pi_0(R)/\mathfrak{m})^{\mathrm{perf}}$, where $\mathfrak{m}$ is a maximal ideal of $\pi_0(R)$ and $(-)^{\mathrm{perf}}$ is the perfection of $\mathbb{F}_p$-algebras. Consequently, we obtain an $\Ee_1$-map 
	\[K(n-1)\otimes(T^n\otimes R)^{hT^n}\longrightarrow K(n-1)\otimes(T^n\otimes(\pi_0(R)/\mathfrak{m})^{\mathrm{perf}})^{hT^n}.\]
	Hence we may assume that $R$ is a perfect field of characteristic $p$. For this, we show that the canonical map $T^n\otimes\mathbb{F}_p\longrightarrow T^n\otimes R$ admits a retraction in $\mathrm{Fun}(BT^n,\mathrm{Sp})$. Since $(T^n\otimes\mathbb{F}_p)^{hT^n}$ has height $n-1$ by \cite[Theorem 1.1]{Zha26}, the desired result follows.
	
	First, \cite[Proposition 2.4]{BSY22} implies that $\mathbb{W}(R)\otimes\mathbb{F}_p\simeq R$. Hence there is an equivalence
	\[T^n\otimes\mathbb{W}(R)\otimes\mathbb{F}_p\simeq T^n\otimes_{\mathbb{F}_p}(\mathbb{W}(R)\otimes\mathbb{F}_p)\simeq T^n\otimes_{\mathbb{F}_p}R\simeq R,\]
	where we write $\otimes_{\mathbb{F}_p}$ for $\otimes_{\mathrm{Mod}_{\mathbb{F}_p}}$. Here the last equivalence follows from the HKR theorem, since $R$ is perfect. Therefore, the natural map $T^n\otimes\mathbb{W}(R)\to\mathbb{W}(R)$ becomes an equivalence after tensoring with $T^n\otimes\mathbb{F}_p$, since $T^n\otimes\mathbb{F}_p$ is a $\mathbb{F}_p$-module and every $\mathbb{F}_p$-module splits as a direct sum of suspensions of $\mathbb{F}_p$. On the other hand, choose an $\mathbb{F}_p$-basis $\{b_i\}_{i\in I}$ of $R$ containing $1$. The Teichm\"uller lifts of $\{b_i\}_{i\in I}$ in $\pi_0(\mathbb{W}(R))$ then give a map 
	\[\bigoplus_{i\in I}\Ss_{(p)}\longrightarrow\mathbb{W}(R).\]
	Using the equivalence $\mathbb{W}(R)\otimes\mathbb{F}_p\simeq R$, we see that this map becomes an equivalence after tensoring with $T^n\otimes\mathbb{F}_p$ by the same reason as above. Since $\THH$ is a symmetric monoidal functor, so is $T^n\otimes-$. Consequently, the above arguments give a zigzag of $T^n$-equivariant equivalences
	\[T^n\otimes R\xleftarrow{\ \simeq \ } (T^n\otimes\mathbb{F}_p)\otimes(T^n\otimes\mathbb{W}(R))\xrightarrow{\ \simeq \ }(T^n\otimes\mathbb{F}_p)\otimes\mathbb{W}(R)\xleftarrow{\ \simeq \ }(T^n\otimes\mathbb{F}_p)\otimes\bigoplus_{i\in I}\Ss_{(p)}.\]
	Since the Teichm\"uller lift of $1$ is $1$, $T^n\otimes R$ contains $T^n\otimes\mathbb{F}_p$ as a direct summand in $\mathrm{Fun}(BT^n,\mathrm{Sp})$.
\end{proof}

We end this section with two applications of Theorem \ref{thm1}. The first application concerns which heights between $m$ and $m+n$ can occur for $(T^n\otimes R)^{hT^n}$, where $R$ is an $\Ee_{\infty}$-ring of height $m$.

\begin{lem}\label{lem:THH and times}
	Let $A$ and $B$ be $\Ee_{\infty}$-rings. For $n\geq1$, there is a $T^n$-equivariant equivalence of $\Ee_{\infty}$-rings
	\[T^n\otimes(A\times B)\simeq (T^n\otimes A)\times(T^n\otimes B).\]
\end{lem}
\begin{proof}
	Applying $\THH$ to the equivalence of symmetric monoidal categories $\mathrm{Perf}(A\times B)\simeq\mathrm{Perf}(A)\times\mathrm{Perf}(B)$, we obtain an equivalence of $\Ee_{\infty}$-rings $T\otimes(A\times B)\simeq (T\otimes A)\times(T\otimes B)$ by the additive theorem; see \cite[Corollary 4.41]{HNS26}. Iterating this $n$ times gives the desired equivalence. 
\end{proof}

\begin{cor}
	For $n\geq1$, $m\geq n-2$ and any $m\leq k\leq m+n$, there is an $\Ee_{\infty}$-ring $R$ of height $m$ such that $(T^n\otimes R)^{hT^n}$ has height $k$.
\end{cor}
\begin{proof}
	First, we consider the case $m\geq1$. When $k=m$ or $m+n$, take $R$ to be $E_m(\Ff)$ or $e_m(\Ff)$, respectively. Then $(T^n\otimes R)^{hT^n}$ has height $m$ or $m+n$ by \cite[Example 5.2]{Zha26} and Theorem \ref{thm1}, respectively.
	
	Now, assume that $n<k<n+m$. Take $R$ to be $E_m(\Ff)\times e_{k-n}(\Ff)$. By Lemma \ref{lem:THH and times}, there is an equivalence $(T^n\otimes R)^{hT^n}\simeq(T^n\otimes E_m(\Ff))^{hT^n}\times(T^n\otimes e_{k-n}(\Ff))^{hT^n}$. Since $k-n>0$, Theorem \ref{thm1} implies that the height of $(T^n\otimes e_{k-n}(\Ff))^{hT^n}$ is $k>n$. Hence we deduce that
	\[\mathrm{height}((T^n\otimes R)^{hT^n})=\max\{\mathrm{height}((T^n\otimes E_m(\Ff))^{hT^n}),\mathrm{height}((T^n\otimes e_{k-n}(\Ff))^{hT^n})\}=k.\]
	Next, assume that $k=n$ or $k=n-1$. Take $R$ to be $E_m(\Ff)\times \Zz$ or $E_m(\Ff)\times \mathbb{F}_p$, respectively. The desired result follows from Theorem 1.1 and Theorem 1.2 in \cite{Zha26}, by the same argument as in the preceding paragraph. 
	
	Finally, we consider the cases $m=-1$ and $m=0$. When $(n,m,k)\neq(1,-1,-1)$, we take
	\[\begin{array}{c|c|c|c|c|c|c}
		\hline
		n & 1 & 1 & 1 & 2 & 2 & 2\\
		m & -1 & 0 & 0 & 0 & 0 & 0\\
		k & 0 & 0 & 1 & 0 & 1 & 2\\
		\hline
		R & \mathbb{F}_p & \mathbb{Q} & \Zz & \mathbb{Q} & \mathbb{Q}\times\mathbb{F}_p & \Zz\\
		\hline
	\end{array}\]
	The same argument as above gives the desired result. When $(n,m,k)=(1,-1,-1)$, we take $R=\mathbb{F}_p^{tC_p}$. By \cite[Theorem I.3.1]{NS18}, there is an $\Ee_{\infty}$-map 
	\[\mathbb{F}_p\longrightarrow\mathbb{F}_p^{hC_p}\longrightarrow\mathbb{F}_p^{tC_p}.\]
	Hence $\mathbb{F}_p^{tC_p}$ is an $\mathbb{F}_p$-module, which implies that $\mathbb{F}_p^{tC_p}$ has height $-1$. We claim that $p[R]$ is zero in $K_0(R)$\footnote{The case $p=2$ follows immediately from \cite[Remark 9.13]{BSY22}.}. Indeed, it is well-known that there is an isomorphism $\mathbb{F}_p[C_p]\cong\mathbb{F}_p[t]/(t^p)$. By \cite[Theorem I.3.1]{NS18}, there is a cofiber sequence of perfect $\mathbb{F}_p^{tC_p}$-modules
	\[(t^{i+1}\mathbb{F}_p[t]/(t^p))^{tC_p}\longrightarrow(t^i\mathbb{F}_p[t]/(t^p))^{tC_p}\longrightarrow \mathbb{F}_p^{tC_p}\]
	for $0\leq i\leq p-1$. Since $t^p\mathbb{F}_p[t]/(t^p)=0$, we deduce that 
	\[p[\mathbb{F}_p^{tC_p}]=[(\mathbb{F}_p[t]/(t^p))^{tC_p}]=[(\mathbb{F}_p[C_p])^{tC_p}]=0.\]
	Since $[\mathbb{F}_p^{tC_p}]$ is the unit of $K_0(\mathbb{F}_p^{tC_p})$, we see that $K(\mathbb{F}_p^{tC_p})\otimes\mathbb{Q}$ is zero. Hence so is $(T\otimes\mathbb{F}_p^{tC_p})^{hT}\otimes\mathbb{Q}$, and the height of $(T\otimes\mathbb{F}_p^{tC_p})^{hT}$ is $-1$. 
\end{proof}

The second application extends Theorem \ref{thm1} from connective $\Ee_{\infty}$-rings to $(-1)$-connective $\Ee_{\infty}$-rings. Again, we need a lemma to prove this.

\begin{lem}\label{lem:connective}
	Let $A$ be a spectrum and $R$ be an $\Ee_{\infty}$-ring equipped with a map of spectra $A\to R$. If both $A$ and $R$ are $(-1)$-connective, then $R\otimes_{\mathrm{Free}_{\Ee_{\infty}}(A)}\Ss$ is also $(-1)$-connective.
\end{lem}
\begin{proof}
	Let $i$ be the inclusion $\mathbb{N}^{\mathrm{ds}}\to \mathbb{N}_<$ and $i_!$ be the left Kan extension $\mathrm{Fun}(\mathbb{N}^{\mathrm{ds}},\mathrm{Sp})\to\mathrm{Fun}(\mathbb{N}_<,\mathrm{Sp})$. Since $i$ is symmetric monoidal, so is $i_!$. Consider $A$ as an object of $\mathrm{Fun}(\mathbb{N}^{\mathrm{ds}},\mathrm{Sp})$ concentrated in weight one. Then we have that
	\begin{equation}\label{eq:fil and gr}
		i_!(\mathbb{F}\mathrm{ree}_{\Ee_{\infty}}(A))_k=\bigoplus_{0\leq i\leq k}(A^{\otimes i})_{h\Sigma_i},\qquad  \mathrm{gr}_k(i_!(\mathbb{F}\mathrm{ree}_{\Ee_{\infty}}(A))\simeq (A^{\otimes k})_{h\Sigma_k},
	\end{equation}
	where $\mathbb{F}\mathrm{ree}_{E_{\infty}}(A)$ is the free $\Ee_{\infty}$-ring generated by $A$ in $\mathrm{Fun}(\mathbb{N}^{\mathrm{ds}},\mathrm{Sp}).$ Write $\mathrm{Free}_{\Ee_{\infty}}(A)$ for the free $\Ee_{\infty}$-ring generated by $A$ in $\mathrm{Sp}$. We have that 
	\[\mathrm{Free}_{\Ee_{\infty}}(A)\simeq\bigoplus_{i\geq0}(A^{\otimes i})_{h\Sigma_i}\simeq\mathrm{colim}_ki_!(\mathbb{F}\mathrm{ree}_{\Ee_{\infty}}(A))_k.\]
	In particular, the $\Ee_{\infty}$-map $\mathrm{Free}_{\Ee_{\infty}}(A)\to R$ adjoint to the given map $A\to R$ corresponds to an $\Ee_{\infty}$-map $i_!(\mathbb{F}\mathrm{ree}_{\Ee_{\infty}}(A))\to \mathrm{const}(R)$ in the category $\mathrm{Fun}(\mathbb{N}_<,\mathrm{Sp})$, where $\mathrm{const}(R)$ is the constant functor with the value $R$. Similarly, the augmentation $\mathrm{Free}_{\Ee_{\infty}}(A)\to \Ss$ corresponds to the $\Ee_{\infty}$-map $i_!(\mathbb{F}\mathrm{ree}_{\Ee_{\infty}}(A))\to \mathrm{const}(\Ss)$. In particular, we can form the following bar construction
	\[B_{\bullet}(\mathrm{const}(R),i_!(\mathbb{F}\mathrm{ree}_{\Ee_{\infty}}(A)),\mathrm{const}(\Ss))\in\mathrm{Fun}(\Delta^{\mathrm{op}},\mathrm{Fun}(\mathbb{N}_<,\mathrm{Sp}))\simeq\mathrm{Fun}(\mathbb{N}_<,\mathrm{Fun}(\Delta^{\mathrm{op}},\mathrm{Sp})).\]
	Since $\mathrm{gr}$ is symmetric monoidal by \cite[Proposition 3.2.1]{Lur15}, Equation \eqref{eq:fil and gr} implies that
	\begin{equation}\label{eq:grbar}
		\mathrm{gr}_*(B_{\bullet}(\mathrm{const}(R),i_!(\mathbb{F}\mathrm{ree}_{\Ee_{\infty}}(A)),\mathrm{const}(\Ss)))\simeq B_{\bullet}(R,\mathrm{gr}_*(i_!(\mathbb{F}\mathrm{ree}_{\Ee_{\infty}}(A)),\Ss),
	\end{equation}
	where we consider $R$ and $\Ss$ as graded spectra concentrated in weight zero. Now, the structure map $d_0$ is induced by the map $R\otimes(A^{\otimes k})_{h\Sigma_k}\to R$, which reduces the filtration by at least one, for $k\geq 1$. Consequently, it follows from the equivalence \eqref{eq:grbar} that
	\begin{equation*}
		\begin{split}
			\mathrm{gr}_*(B_{\bullet}(\mathrm{const}(R),i_!(\mathbb{F}\mathrm{ree}_{\Ee_{\infty}}(A)),\mathrm{const}(\Ss)) & \simeq R\otimes B_{\bullet}(\Ss,\mathrm{gr}_*(i_!(\mathbb{F}\mathrm{ree}_{\Ee_{\infty}}(A)),\Ss)\\
			& \simeq R\otimes B_{\bullet}(\Ss,\mathbb{F}\mathrm{ree}_{\Ee_{\infty}}(A),\Ss)\\
		\end{split}
	\end{equation*}
	since $\mathrm{gr}_*\circ i_!$ is equivalent to the identity. Taking the geometric realization implies that 
	\[\mathrm{gr}_*(\mathrm{const}(R)\otimes_{i_!(\mathbb{F}\mathrm{ree}_{\Ee_{\infty}}(A))}\mathrm{const}(\Ss))\simeq R\otimes\mathbb{F}\mathrm{ree}_{\Ee_{\infty}}(\Sigma A),\]
	since $\mathbb{F}\mathrm{ree}_{\Ee_{\infty}}(-)$ is a left adjoint. Hence we see that $\mathrm{gr}_k(\mathrm{const}(R)\otimes_{i_!(\mathbb{F}\mathrm{ree}_{\Ee_{\infty}}(A))}\mathrm{const}(\Ss))$ is $(-1)$-connective, which implies that $(\mathrm{const}(R)\otimes_{i_!(\mathbb{F}\mathrm{ree}_{\Ee_{\infty}}(A)}\mathrm{const}(\Ss))_k$ is $(-1)$-connective for each $k\geq0$. Since $\mathbb{N}_{<}$ is sifted, $R\otimes_{\mathrm{Free}_{\Ee_{\infty}(A)}}\Ss\simeq\mathrm{colim}_k(\mathrm{const}(R)\otimes_{i_!(\mathbb{F}\mathrm{ree}_{\Ee_{\infty}}(A))}\mathrm{const}(\Ss))$. This implies the desired result, since taking colimit preserves connectivity.
\end{proof}

\begin{cor}
	Let $R$ be a $(-1)$-connective $\Ee_{\infty}$-ring of height $m\geq1$. For any $n\geq1$, the height of $(T^n\otimes R)^{hT^n}$ is $m+n$.
\end{cor}
\begin{proof}
	Since $R$ has height $m$, $L_{T(m)}R$ is nonzero. Hence we obtain an $\Ee_{\infty}$-map $R\to E_m(L)$ for some algebraically closed field $L$, by \cite[Corollary 5.2]{BSY22}. Using the small object argument, we show that there is an $\Ee_{\infty}$-ring $R^{\prime}$ of height $m$ such that the above map factors as 
	\[R\longrightarrow R^{\prime}\longrightarrow E_m(L).\]
	Take $R_0$ to be $R$. Suppose that we have constructed an $\Ee_{\infty}$-ring $R_i$ and a factorization $R\to R_i\to E_m(L)$ such that $R_i$ is $(-1)$-connective. Choose a representative $x\colon\Sigma^{-1}\Ss\to R_i$ for each element $x$ of $\pi_{-1}(R_i)$. Since $R_i$ is an $\Ee_{\infty}$-ring, $x$ is adjoint to an $\Ee_{\infty}$-map $\mathrm{Free}_{\Ee_{\infty}}(\Sigma^{-1}\Ss)\to R_i$, which is still denoted by $x$. Since $\pi_{-1}(E_m(L))=0$, the composition 
	\[\mathrm{Free}_{\Ee_{\infty}}(\Sigma^{-1}\Ss)\xrightarrow{\ x\ } R_i\longrightarrow E_m(L)\]
	factors as
	\[\mathrm{Free}_{\Ee_{\infty}}(\Sigma^{-1}\Ss)\longrightarrow\Ss\longrightarrow E_m(L).\]
	Set $A_i=\bigotimes_{x\in\pi_{-1}(R_i)}\mathrm{Free}_{\Ee_{\infty}}(\Sigma^{-1}\Ss)$ and $R_{i+1}=R_i\otimes_{A_i}\Ss$. By Lemma \ref{lem:connective}, $R_{i+1}$ is $(-1)$-connective. Moreover, the composition of $x$ and the canonical map $R_i\to R_{i+1}$ is zero. We take $R^{\prime}$ to be $\mathrm{colim}_iR_i$. Since $\Sigma^{k}\Ss$ is compact for any $k\in\Zz$, any map $\Sigma^k\Ss\to R^{\prime}$ factors through $R_i$ for some $i$. When $k<-1$, this map is zero, since $R_i$ is $(-1)$-connective. When $k=-1$, this map is zero after composing with the canonical map $R_i\to R_{i+1}$ by construction. Hence $R^{\prime}$ is connective.
	
	Since $R$, $R^{\prime}$ and $E_m(L)$ are $\Ee_{\infty}$-rings, \cite[Lemma 2.3]{LMMT24} implies that 
	\[m=\mathrm{height}(E_m(L))\leq\mathrm{height}(R^{\prime})\leq\mathrm{height}(R)=m.\]
	Hence the height of $R^{\prime}$ is also $m$. Since $R^{\prime}$ is connective, Theorem \ref{thm1} implies that the height of $(T^n\otimes R^{\prime})^{hT^n}$ is $m+n$, which further implies that $(T^n\otimes R)^{hT^n}$ has height at least $m+n$. On the other hand, there is an $\Ee_{\infty}$-map $\tau_{\geq0}R\to R$. Since $\tau_{<0}R$ is bounded above, $K(k)\otimes\tau_{<0}R$ is zero for $k\geq1$, by Lemma 2.2 and Lemma 2.3 in \cite{LMMT24}. Hence $\tau_{\geq0}R$ also has height $m$. Applying $(T^n\otimes-)^{hT^n}$ to this $\Ee_{\infty}$-map and using Theorem \ref{thm1} again show that $(T^n\otimes R)^{hT^n}$ has height at most $m+n$.
\end{proof}

\section{The proof of Theorem \ref{thm2}}\label{sec5}
As noted in the introduction, the case $m=0$ behaves differently. Therefore, we treat it separately in this section. The idea is similar to that used in the proof of Theorem \ref{thm1}. More precisely, for a $p$-local connective $\Ee_{\infty}$-ring $R$ of height zero with $\pi_0(R)_p^{\wedge}[1/p]$ nonzero, we construct an $\Ee_{\infty}$-ring $B$, together with an $\Ee_{\infty}$-map $R\to B$, such that the height of $(T^n\otimes B)^{hT^n}$ is easier to compute. For this, we need the following lemma.

\begin{lem}\label{lem:A to B}
	Let $A$ be a commutative ring such that $A_p^{\wedge}[1/p]$ is nonzero. Then there is a cospan of commutative rings
	\[A\longrightarrow B\xleftarrow{\ \ \epsilon \ \ } W(B^{\flat}),\]
	where $(-)^{\flat}$ is computed by the inverse limit along Frobenius on $B/p$, and $\epsilon$ is induced by the counit of the adjunction in \cite[Proposition 2.2]{BSY22}. Moreover, $\epsilon$ is surjective with kernel generated by a nonzerodivisor $\xi$.
\end{lem}
\begin{proof}
	Let $I\subseteq A_p^{\wedge}$ be the contraction of a maximal ideal of $A_p^{\wedge}[1/p]$. Then $A_p^{\wedge}/I$ is a domain with $p\neq0$ and $(A_p^{\wedge}/I)[1/p]$ is a field. Hence, for any nonzero $a\in A_p^{\wedge}/I$, we have that $a$ is nonzero in $(A_p^{\wedge}/I)[1/p]$ and $a^{-1}=b/p^k$ in $(A_p^{\wedge}/I)[1/p]$ for some $k\geq0$ and $b\in A_p^{\wedge}/I$. In particular, $ab=p^k$ in $A_p^{\wedge}/I$. Moreover, since $1-pa$ is invertible for any $a\in A_p^{\wedge}/I$, $p$ lies in the Jacobson radical of $A_p^{\wedge}/I$. 
	
	Now, using \cite[Lemma 10.50.2]{Stacks}, we may choose a maximal ideal $\mathfrak{m}$ of $A_p^{\wedge}/I$ and a valuation ring of $\mathrm{Frac}(A_p^{\wedge}/I)$ dominating $(A_p^{\wedge}/I)_{\mathfrak{m}}$. Denote the associated valuation by $v\colon\mathrm{Frac}(A_p^{\wedge}/I)\to\Gamma\cup\{\infty\}$, where $\Gamma$ is a totally ordered abelian group. Since $p$ lies in the Jacobson radical of $A_p^{\wedge}/I$, $v(p)$ is greater than zero. Moreover, write $\Gamma$ multiplicatively with the order reversed, and write $|-|_v$ for the corresponding multiplicative valuation. For nonzero elements $a,c\in A_p^{\wedge}/I$, the above argument gives $k\geq0$ and $b\in A_p^{\wedge}/I$ such that $ab=p^k$. Hence, for $n>k$, we have that
	\[|p|_v^n<|p|_v^k=|a|_v|b|_v\leq |a|_v\leq |a/c|_v,\]
	since $v(b),v(c)\geq v(1)=0$ and $v(p)>v(1)=0$. By \cite[Lemma 4.2.2]{SW20}, there is an order preserving group homomorphism $\rho\colon\Gamma\to\mathbb{R}_{>0}$ which sends $|p|_v$ to $1/2$. Hence we obtain an $\mathbb{R}$-valued valuation given by $w(x):=-\log_2\bigl(\rho(|x|_v)\bigr)$ for $x\in\mathrm{Frac}(A_p^{\wedge}/I)^{\times}$. Set $|x|:=p^{-w(x)}$ for $x\in\mathrm{Frac}(A_p^{\wedge}/I)^{\times}$ and $|0|=0$. Then $|xy|=|x||y|$ and $|x+y|\leq\max\{|x|,|y|\}$. That is, $|-|$ defines a nonarchimedean absolute value on $\mathrm{Frac}(A_p^{\wedge}/I)$ with $|p|=p^{-1}$. Let $C$ be the completion of $\mathrm{Frac}(A_p^{\wedge}/I)$ with respect to $|-|$. By \cite[Appendix A, Theorem 3]{Bos14}, $|-|$ extends uniquely to an algebraic closure $\overline{C}$ of $C$. Let $D$ be the completion of $\overline{C}$ with respect to $|-|$, and set $B=\{x\in D\ |\ |x|\leq1\}$. Then $D$ is algebraically closed by \cite[Appendix A, Lemma 6]{Bos14}, and, as the closed unit ball in $D$, $B$ is complete. For a nonzero element $a\in A_p^{\wedge}/I$, since the embedding $\mathrm{Frac}(A_p^{\wedge}/I)\hookrightarrow D$ preserves this absolute value and $w(a)\geq0$, we have that $|a|\leq1$. Hence $B$ contains the image of $A$, which gives the desired map $A\to B$. 
	
	Finally, note that, since $|p|=p^{-1}$, the ideals $p^kB=\{x\in B\ |\ |x|\le p^{-k}\}$ form a neighborhood basis of zero. Hence $B$ is $p$-adically complete, which implies that $B$ is $p$-complete as a spectrum. Consequently, the adjunction in \cite[Proposition 2.2]{BSY22} induces a map $\epsilon\colon W(B^{\flat})\to B$. Since $D$ is algebraically closed, there is an element $\pi$ with $\pi^p=p$ and $|\pi|=|p|^{1/p}<1$. In particular, $\pi\in B$ is a nonzerodivisor, and $\pi^kB=\{x\in B\ |\ |x|\leq |\pi|^k\}$, which form a neighborhood basis of zero. Hence $B$ is $\pi$-adically complete. Moreover, the map
	\[B/\pi\xrightarrow{\ x\mapsto x^p} B/p\]
	is surjective, since $D$ is algebraically closed and $|x|^p\leq1$ implies that $|x|\leq1$. It is also injective, since $x^p\in pB$ if and only if $|x|\leq |p|^{1/p}=|\pi|$. Consequently, Lemma 3.9 and Lemma 3.10 in \cite{BMS18} imply that $\epsilon$ is surjective with kernel generated by a nonzerodivisor $\xi$. 
\end{proof}

\begin{thm}
	Let $R$ be a $p$-local connective $\Ee_{\infty}$-ring of height zero. For any $n\geq1$, we have that
	\begin{equation*}
		\mathrm{height}((T^n\otimes R)^{hT^n})=
		\begin{cases}
			0, & p\in\pi_0(R)^{\times}, \\
			n-1, & \text{$p\notin\pi_0(R)^{\times}$ and $\Rr[1/p]=0$,} \\
			n, & \Rr[1/p]\neq0.
		\end{cases}
	\end{equation*}
\end{thm}
\begin{proof}
	If $p$ is invertible in $\pi_0(R)$, then $R$ is rational. Hence the height of $(T^n\otimes R)^{hT^n}$ is zero as discussed in the introduction.
	
	Therefore, we assume that $p$ is not invertible in $\pi_0(R)$. First, we consider the case $\Rr[1/p]=0$. Note that $\Rr[1/p]$ is zero if and only if $p^N\pi_0(R)=p^{N+1}\pi_0(R)$ for some $N>0$. Indeed, if $\Rr[1/p]$ is zero, then $p^N\cdot1=0$ in $\Rr$ for some $N>0$. Hence $p^N\cdot\overline{1}=0$ in $\pi_0(R)/p^{N+1}\pi_0(R)$, which implies that $p^N\pi_0(R)\subseteq p^{N+1}\pi_0(R)$. Conversely, if $p^N\pi_0(R)=p^{N+1}\pi_0(R)$, then $p^N(1-px)=0$ for some $x\in\pi_0(R)$. Since $1-px$ is invertible in $\Rr$, we deduce that $p^N=0$ in $\Rr$. Moreover, for $x$ and $N$ above, we have that
	\[1=(1-px)(1+px+\cdots+(px)^{N-1})+(px)^N.\]
	That is, the ideals $(p^N)$ and $(1-px)$ are coprime. Since $(p^N)(1-px)=0$ by the discussion above, we obtain the isomorphism 
	\[\pi_0(R)\cong\pi_0(R)/((p^N)\cap(1-px))\cong\pi_0(R)/p^N\times \pi_0(R)/(1-px).\]
	Let $e\in\pi_0(R)$ be the idempotent corresponding to the summand $\pi_0(R)/p^N$. Then there is an equivalence 
	\[R\simeq R[e^{-1}]\times R[(1-e)^{-1}].\]
	Since $p$ is not invertible in $\pi_0(R)$, $\pi_0(R)/p^N$ is nonzero, which implies that $R[e^{-1}]$ has height $-1$. Since $p$ is invertible in $\pi_0(R[(1-e)^{-1}])\cong\pi_0(R)/(1-px)$ and $R$ has height zero, $R[(1-e)^{-1}]$ is also nonzero and rational. Hence Lemma \ref{lem:THH and times} implies that
	\[\mathrm{height}((T^n\otimes R)^{hT^n})=\max\{\mathrm{height}((T^n\otimes R[e^{-1}])^{hT^n}),\mathrm{height}((T^n\otimes R[(1-e)^{-1}])^{hT^n})\}=n-1\]
	by Theorem \ref{thm1} and the discussion in the introduction.
	
	Consequently, it remains to consider the case $\Rr[1/p]\neq0$. As observed in the introduction, we need to show that $K(n)\otimes(T^n\otimes R)^{hT^n}$ is nonzero. By Lemma \ref{lem:A to B}, there exist a commutative ring $B$ and an $\Ee_1$-map
	\[K(n)\otimes(T^n\otimes R)^{hT^n}\longrightarrow K(n)\otimes(T^n\otimes B)^{hT^n}.\]
	Hence we only need to show that $K(n)\otimes(T^n\otimes B)^{hT^n}$ is nonzero. Set $U=\MU_p\otimes_{\Ss_p}\mathbb{W}(B^{\flat})$. By the same computation as in Lemma \ref{lem:pi(W0)}, we have that
	\[\pi_*(U)\cong W(B^{\flat})[v_i,x_j\ |\ j\neq p^k-1].\]
	In particular, $U$ is connective with $\pi_0(U)\cong W(B^{\flat})$. By Lemma \ref{lem:A to B}, we obtain an $\Ee_{\infty}$-map $U\to B$, which is surjective on homotopy groups with kernel generated by the regular sequence $\{\xi,v_i,x_j\}$, where $\xi$ is the nonzerodivisor given in Lemma \ref{lem:A to B}. Moreover, we have a commutative diagram
	\[\xymatrix{
	\MU \ar[d] \ar[r] & \Zz_{(p)} \ar[d] \ar[r] & \Zz_{(p)} \ar[d]\\
	U \ar[r] & W(B^{\flat}) \ar[r] & B,}\]
	which induces an $\Ee_{\infty}$-map $T^n\otimes_{\MU}\Zz_{(p)}\to T^n\otimes_UB$. Hence the relation in \cite[Lemma 2.3]{Zha26} also holds in $\pi_*((T^n\otimes_UB)/p)$. Now, the same argument as in the proof of \cite[Theorem 1.2]{Zha26} with $\MU$ and $\Zz$ replaced by $U$ and $B$ implies that $K(n)\otimes(T^n\otimes B)^{hT^n}$ is nonzero.
\end{proof}

\section{The supports of $(\int_{T^2}\BP{n})^{hT^2}$ and $(\int_{T^3}\BP{n})^{hT^3}$}\label{sec6}
In this section, we compute the supports of $(\int_{T^2}\BP{n})^{hT^2}$ and $(\int_{T^3}\BP{n})^{hT^3}$. First, we use Theorem \ref{thm1} to show that these two supports contain $\{0,\cdots,n+2\}$ and $\{0,\cdots,n+3\}$, respectively. The following proposition is similar to \cite[Construction 5.10]{BHLS23}, but we do not require the $\Ee_3$-$\MU$-form on $\BP{n}$ to be the one constructed in \cite[\S2]{HW22}. 

\begin{prop}\label{prop:bptoe}
	There is an $\Ee_3$-map $\BP{n}\to e_n(\Ff)$.
\end{prop}
\begin{proof}
	Since $p^n-1$ is invertible in the $\Ee_{\infty}$-ring $\MU[v_n^{-1}]$, \cite[Example 31]{Law20} implies that there is  an $\Ee_{\infty}$-$\MU$-algebra $\MU[v_n^{-1}][x]/[x^{p^n-1}-v_n]$. Base changes along its unit map and $\Ss\to\mathbb{W}(\Ff)$, followed by completion with respect to $(p,v_1,\cdots,v_{n-1})$, give an $\Ee_3$-ring 
	\[(\mathbb{W}(\Ff)\otimes\BP{n}\otimes_{\MU}\MU[v_n^{-1}][x]/[x^{p^n-1}-v_n])_{(p,v_1,\cdots,v_{n-1})}^{\wedge}.\]
	The same argument as in the proof of Lemma \ref{lem:pi(W0)} implies that its homotopy group is isomorphic to $W(\Ff)[[u_1,\cdots,u_{n-1}]][u^{\pm}]$ with $u_i=v_iu^{-(p^i-1)}$.
	
	By Lemma \ref{lem:e3}, it remains to identify the formal group law $F$ over this $\Ee_3$-ring. By \cite[Proposition 1]{Lur-Lecture13}, the coefficient of $x^{p^i}$ in $[p]_{F}(x)$ is $b_i=a_iu_i+F_i(u_1,\cdots,u_{i-1})$ with $a_i\in\Zz_{(p)}^{\times}$, for $1\leq i\leq n-1$. Moreover, $b_n=a_n+F_n(u_1,\cdots,u_{n-1})$ with $a_n\in\Zz_{(p)}^{\times}$. Here $F_i$ is a polynomial with zero constant term. Hence, for $1\leq i\leq n-1$, we obtain that $u_i=a_i^{-1}(b_i-F_i(u_1,\cdots,u_{i-1}))$, which implies that 
	\[W(\Ff)[[u_1,\cdots,u_{n-1}]][u^{\pm}]\cong W(\Ff)[[b_1,\cdots,b_{n-1}]][u^{\pm}].\]
	Since $b_n$ is a unit in $W(\Ff)[[b_1,\cdots,b_{n-1}]][u^{\pm}]$, $F\otimes_{W(\Ff)[[u_1,\cdots,u_{n-1}]]}\Ff$ has height $n$, the result follows from \cite[Theorem 5]{Lur-Lecture21}.
\end{proof}

\begin{cor}\label{cor:heightbp}
	For $m=2$ or $3$, $K(k)\otimes(\int_{T^m}\BP{n})^{hT^m}$ is nonzero for $0\leq k\leq m+n$. In particular, $T(k)\otimes(\int_{T^m}\BP{n})^{hT^m}$ is also nonzero for $0\leq k\leq m+n$.
\end{cor}
\begin{proof}
	By Proposition \ref{prop:bptoe}, tensoring with $K(k)$ gives an $\Ee_0$-map 
	\[K(k)\otimes(\int_{T^m}\BP{n})^{hT^m}\to K(k)\otimes(\int_{T^m}e_n(\Ff))^{hT^m}.\]
	Since $(\int_{T^m}e_n(\Ff))^{hT^m}$ is an $\Ee_{\infty}$-ring of height $m+n$ by Theorem \ref{thm1}, the unit map $\Ss\to K(k)\otimes(\int_{T^m}e_n(\Ff))^{hT^m}$ is nonzero for $0\leq k\leq m+n$. Hence $\Ss\to K(k)\otimes(\int_{T^m}\BP{n})^{hT^m}$ is also nonzero for $0\leq k\leq m+n$. In particular, $K(k)\otimes(\int_{T^m}\BP{n})^{hT^m}$ is nonzero for $0\leq k\leq m+n$.
\end{proof}

Next, we prove that the supports of $(\int_{T^2}\BP{n})^{hT^2}$ and $(\int_{T^3}\BP{n})^{hT^3}$ are exactly $\{0,\cdots,n+2\}$ and $\{0,\cdots,n+3\}$, respectively. We first consider the case $m=2$. The case $m=3$ is more difficult, since $(\int_{T^3}\BP{n})^{hT^3}$ may not admit an $\Ee_1$-ring structure. Since $\BP{n}$ is an $\Ee_3$-ring, $(\int_{T^2}\BP{n})^{hT^2}$ admits an $\Ee_1$-ring structure. By \cite[Theorem 7.4]{BGT14}, the cyclotomic trace map, the canonical natural transformation $\mathrm{TC}(-)\to\THH(-)^{hT}$, and the comparison map of limits give an $\Ee_1$-map
\[K(k)\otimes K^{(2)}(\BP{n})\to K(k)\otimes(\int_{T^2}\BP{n})^{hT^2}.\]
Lemma 2.3 and Corollary 4.1 in \cite{LMMT24} then imply that $K(k)\otimes  K^{(2)}(\BP{n})\simeq0$ for $k>n+2$. Consequently, $K(k)\otimes(\int_{T^2}\BP{n})^{hT^2}$ is also zero for $k>n+2$. Since $(\int_{T^2}\BP{n})^{hT^2}$ is an $\Ee_1$-ring, \cite[Lemma 2.3]{LMMT24} implies that $T(k)\otimes(\int_{T^2}\BP{n})^{hT^2}$ is also zero for $k>n+2$.

Therefore, to prove Theorem \ref{thm3}, it remains to prove the following proposition.

\begin{prop}
	For $k>n+3$, $T(k)\otimes(\int_{T^3}\BP{n})^{hT^3}\simeq0$.
\end{prop}
\begin{proof}
	By \cite[Corollary 5.3.1.15]{HA}, there exists a center $\mathfrak{Z}_{\Ee_3}(\BP{n})\in\mathrm{Alg}_{\Ee_4}(\mathrm{Sp})$ such that $\BP{n}$ is a $\mathfrak{Z}_{\Ee_3}(\BP{n})$-module in the category $\mathrm{Alg}_{\Ee_3}(\mathrm{Sp})$. In particular, the composition of $\Ee_3$-maps
	\begin{equation}\label{eq:center}
		\BP{n}\xrightarrow{\eta\otimes
	id_{\BP{n}}} \mathfrak{Z}_{\Ee_3}(\BP{n})\otimes\BP{n}\longrightarrow\BP{n}
	\end{equation}
	is the identity map of $\BP{n}$. 
	
	We claim that $T(k)\otimes \mathfrak{Z}_{\Ee_3}(\BP{n})\simeq0$ for $k>n$. By Remark 5.3.1.32 and Theorem 7.1.2.1 in \cite{HA} and \cite[Remark 2.2.3]{HW22}, we have that 
	\[\mathfrak{Z}_{\Ee_3}(\BP{n})\simeq \mathrm{map}_{\int_{S^2}\BP{n}}(\BP{n},\BP{n}).\]
	By \cite[Proposition 6.2.1]{HW22}, there is a finite spectrum $M$ such that 
	\begin{equation}\label{eq:MBP}
		M\otimes\BP{n}\simeq\bigoplus_{i=1}^m\Sigma^{d_i}\mathbb{F}_p
	\end{equation}
	for some $m>1$. Hence we obtain that 
	\begin{equation*}
		\begin{split}
			& M\otimes\mathfrak{Z}_{\Ee_3}(\BP{n})\simeq M\otimes\mathrm{map}_{\int_{S^2}\BP{n}}(\BP{n},\BP{n})\simeq\mathrm{map}_{\int_{S^2}\BP{n}}(\BP{n},M\otimes\BP{n})\\
			& \simeq\lim_{\Delta}\left[\mathrm{map}(\BP{n},M\otimes\BP{n})\rightrightarrows\mathrm{map}(\BP{n}\otimes\int_{S^2}\BP{n},M\otimes\BP{n})\threerightarrows\cdots\right]
		\end{split}
	\end{equation*}
	Since $\int_{S^2}\BP{n}$ is connective, we have that $\pi_k(\mathrm{map}(\BP{n}\otimes(\int_{S^2}\BP{n})^{\otimes l},M\otimes\BP{n}))=0$ for $l\geq0$ and $k>\max_i(d_i)$, by the equivalence \eqref{eq:MBP}. That is, $M\otimes\mathfrak{Z}_{\Ee_3}(\BP{n})$ is bounded above, which implies that $T(k)\otimes M\otimes\mathfrak{Z}_{\Ee_3}(\BP{n})\simeq0$ for $k>0$. By \cite[Proposition 6.2.1]{HW22}, $K(k)\otimes M$ is nonzero for $k>n$. Hence, for $k>n$, $K(k)\otimes T(k)\otimes M\otimes\mathfrak{Z}_{\Ee_3}(\BP{n})\simeq0$ splits as a direct sum of suspensions of $K(k)\otimes\mathfrak{Z}_{\Ee_3}(\BP{n})$, which implies that $K(k)\otimes\mathfrak{Z}_{\Ee_3}(\BP{n})\simeq0$. Since $\mathfrak{Z}_{\Ee_3}(\BP{n})$ is an $\Ee_4$-ring, we see that $T(k)\otimes \mathfrak{Z}_{\Ee_3}(\BP{n})\simeq0$ for $k>n$, by \cite[Lemma 2.3]{LMMT24}.
	
	Now, by \cite[Theorem 7.4]{BGT14}, the cyclotomic trace map, the canonical natural transformation $\mathrm{TC}(-)\to\THH(-)^{hT}$, and the comparison map of limits give an $\Ee_1$-map 
	\[K^{(3)}(\mathfrak{Z}_{\Ee_3}(\BP{n}))\to(\int_{T^3}\mathfrak{Z}_{\Ee_3}(\BP{n}))^{hT^3}.\]
	By the above claim, using \cite[Corollary 4.1]{LMMT24} three times, we obtain that $T(k)\otimes K^{(3)}(\mathfrak{Z}_{\Ee_3}(\BP{n}))\simeq0$ for $k>n+3$. Using \cite[Lemma 2.3]{LMMT24} again, we have that $K(k)\otimes K^{(3)}(\mathfrak{Z}_{\Ee_3}(\BP{n}))\simeq0$ for $k>n+3$, which implies that $K(k)\otimes(\int_{T^3}\mathfrak{Z}_{\Ee_3}(\BP{n}))^{hT^3}\simeq0$ for $k>n+3$. Hence so is $T(k)\otimes(\int_{T^3}\mathfrak{Z}_{\Ee_3}(\BP{n}))^{hT^3}$, since $(\int_{T^3}\mathfrak{Z}_{\Ee_3}(\BP{n}))^{hT^3}$ is an $\Ee_1$-ring. Finally, applying $(\int_{T^3}-)^{hT^3}$ to the map \eqref{eq:center}, we see that 
	\[(\int_{T^3}\BP{n})^{hT^3}\longrightarrow(\int_{T^3}\mathfrak{Z}_{\Ee_3}(\BP{n}))^{hT^3}\otimes(\int_{T^3}\BP{n})^{hT^3}\longrightarrow(\int_{T^3}\BP{n})^{hT^3}\]
	is the identity map. That is, $T(k)\otimes(\int_{T^3}\BP{n})^{hT^3}$ is a retract of $T(k)\otimes(\int_{T^3}\mathfrak{Z}_{\Ee_3}(\BP{n}))^{hT^3}\otimes(\int_{T^3}\BP{n})^{hT^3}$, which gives the desired result.
\end{proof}

\begin{cor}\label{cor:supp(K(BPn))}
	For $m=2,3$, the support of $K^{(m)}(\BP{n})$ is $\{0,1,\cdots,m+n\}$.
\end{cor}
\begin{proof}
	Again, by \cite[Theorem 7.4]{BGT14}, the cyclotomic trace map, the canonical natural transformation $\mathrm{TC}(-)\to\THH(-)^{hT}$, and the comparison map of limits give an $\Ee_0$-map
	\[K(k)\otimes K^{(m)}(\BP{n})\longrightarrow K(k)\otimes(\int_{T^m}\BP{n})^{hT^m}.\]
	In other words, the composition
	\[\Ss\to K(k)\otimes K^{(m)}(\BP{n})\to K(k)\otimes(\int_{T^m}\BP{n})^{hT^m}\]
	is the unit map. By Corollary \ref{cor:heightbp} and its proof, this composition is nonzero for $0\leq k\leq m+n$. Hence $K(k)\otimes K^{(m)}(\BP{n})$ is also nonzero for $0\leq k\leq m+n$. In particular, so is $T(k)\otimes K^{(m)}(\BP{n})$. Since \cite[Corollary 4.1]{LMMT24} implies that $T(k)\otimes K^{(m)}(\BP{n})$ is zero for $k>m+n$, the support of $K^{(m)}(\BP{n})$ is exactly $\{0,1,\cdots,m+n\}$.
\end{proof}

\end{document}